\documentclass[a4paper,10pt]{amsart}
\usepackage[foot]{amsaddr}

\usepackage[usenames,dvipsnames]{xcolor}
\usepackage{
 amssymb
,amsfonts
,amsmath
,url
,pigpen
,enumitem
,stmaryrd
,pifont
,hyphenat
,epigraph
,tikz
,etoolbox
,mparhack
,mathtools
,microtype
,lmodern
,datetime
,cite
,CJKutf8
,kodi
,marvosym
,chessfss
}

\usetikzlibrary{fit}

\usepackage[all,2cell,cmtip]{xy}\UseAllTwocells

\renewcommand{\textbf}[1]{\text{\fontseries{b}\selectfont{\upshape #1}}}

\usepackage{xspace}
  \providecommand{\abbrv}[1]{#1.\@\xspace}
  \providecommand{\ie}       {\abbrv{i.e}}

  \providecommand{\eg}       {\abbrv{e.g}}
  \providecommand{\adef}     {\abbrv{Def}}
  \providecommand{\acor}     {\abbrv{Cor}}
  \providecommand{\aprop}    {\abbrv{Prop}}
  \providecommand{\athm}     {\abbrv{Thm}}

\usepackage{cjhebrew}
\usepackage[spanish, greek, english]{babel}

\usepackage[utf8]{inputenc}
\usepackage[T1]{fontenc}

\newtheoremstyle{reference}
   {}                
   {}                
   {}              
   {}                      
   {\fontseries{b}\selectfont}              
   {:}                     
   {.2em}                  
   {\thmname{#1}           
    \thmnumber{#2}         
    \thmnote{{\sc [#3]}}} 

\theoremstyle{reference}
  \newtheorem{theorem}{Theorem}[section]
  \newtheorem{lemma}[theorem]{Lemma}
  \newtheorem{proposition}[theorem]{Proposition}
  \newtheorem{example}[theorem]{Example}
  \newtheorem{remark}[theorem]{Remark}
  \newtheorem{definition}[theorem]{Definition}
  \newtheorem{corollary}[theorem]{Corollary}
  \newtheorem{notat}[theorem]{Notation}

  \newtheorem*{theorem*}{Theorem}
  \newtheorem*{lemma*}{Lemma}
  \newtheorem*{proposition*}{Proposition}
  \newtheorem*{example*}{Example}
  \newtheorem*{question*}{Exercise}
  \newtheorem*{remark*}{Remark}
  \newtheorem*{definition*}{Definition}
  \newtheorem*{corollary*}{Corollary}
  \newtheorem*{notat*}{Notation}
  \newtheorem*{scholium*}{Scholium}

\usepackage{hyperref}
\hypersetup{
  pdftitle={hearts and towers},
  pdfauthor={D. Fiorenza and F. Loregian and G. Marchetti},
  colorlinks=true,
  linkcolor=blue!40!black,
  citecolor=blue!40!black
}

\makeatletter
  \def\@cite#1#2{[\textbf{#1}\if@tempswa , #2\fi]}
  \def\@biblabel#1{[\textsf{#1}]}
\makeatother

\newcommand{\refbf}[1]{\textbf{\ref{#1}}}

\newcommand{\orth}{\boxslash}

\renewcommand{\bowtie}{%
  \mathchoice%
  {\lhd\hskip-1.8pt\rhd}%
  {\lhd\hskip-1.8pt\rhd}%
  {\lhd\hskip-1.8pt\rhd}%
  {\lhd\hskip-1.8pt\rhd}%
}

\usepackage[titletoc]{appendix}
\usepackage[bottom=6cm,
			left=4cm,
			right=4.5cm,
			top=4.5cm]{geometry}

\newcommand{\xto}[1]{\xrightarrow{#1}}

\renewcommand{\to}{\longrightarrow}

\newcommand{\tee}{\mathfrak{t}}
\newcommand{\po}{\ar@{}[dr]|(.7){\text{ }}} 
\newcommand{\pb}{\ar@{}[dr]|(.3){\text{\pigpenfont J}}}
\newcommand{\pp}{\ar@{}[dr]|{\text{ }}} 
\newcommand{\heart}{\text{\Heart}}

\newcommand{\var}[2]{\left[\begin{smallmatrix} #1 \\ \downarrow \\ #2 \end{smallmatrix}\right]}
\newcommand{\varnobkt}[2]{\begin{smallmatrix} #1 \\ \downarrow \\ #2 \end{smallmatrix}}

\newcommand{\cate}[1]{\EuScript{#1}}
\DeclareMathAlphabet\EuScript{U}{eus}{m}{n}
\SetMathAlphabet\EuScript{bold}{U}{eus}{b}{n}

\newcommand{\cA}  {\cate{A}}
\newcommand{\cB}  {\cate{B}}
\newcommand{\C}   {\EuScript{C}}

\newcommand{\A}   {\mathcal{A}}
\newcommand{\E}   {\mathcal{E}}
\newcommand{\M}		{\mathcal{M}}
\renewcommand{\O} {\mathcal{O}}
\newcommand{\fF}  {\mathbb{F}}
\newcommand{\R}   {\mathbb{R}}

\newcommand{\Z}   {\mathbb{Z}}

\newcommand{\T}   {\mathcal{T}}

\renewcommand{\phi}{\varphi}

\newlength{\seplen}
\newlength{\sepwid}
\def\firstblank{\,\rule{\seplen}{\sepwid}\,}
\def\secondblank{\firstblank\llap{\raisebox{2pt}{\firstblank}}}

\DeclareMathOperator{\im}{\textsf{im}}
\DeclareMathOperator{\coim}{\textsf{coim}} 			
\DeclareMathOperator{\coker}{coker}
\DeclareMathOperator{\fib}{fib}
\DeclareMathOperator{\cofib}{cofib}

\newcommand{\smallcap}[1]{\text{\scshape #1}}
  \providecommand{\fs}{\smallcap{fs}}
  \providecommand{\ts}{\smallcap{ts}}  
  \providecommand{\pf}{\smallcap{pf}}
  \providecommand{\tot}{\smallcap{tt}}

\newcommand{\EmailTo}[1]{\href{mailto:#1}{\sf #1}}

\definecolor{semilightgray}{rgb}{0.65, 0.65, 0.65}

\usepackage{turnstile}

\def\lu{{(L,U)}}

\newcommand{\ordered}[1]{[\mathbf{#1}]}
  \def\uno{\ordered{1}}
  \def\zero{\textbf{0}}
  \def\ordk{\ordered{k}}

\newcommand{\Pos}{\underline{\cate{P}\text{os}}}
\newcommand{\Tos}{\underline{\cate{T}\text{os}}}
\def\ZPos{\Z\text{-}\Pos}
\def\Sp{\cate{S}\text{p}}
\def\weaves{\text{\Denarius}}
\newcommand{\japanese}[2][bsmi]{\text{\begin{CJK*}{UTF8}{#1}#2\end{CJK*}}}
\newcommand{\slicings}{\japanese{切}}

\newlength\DefaultEpiWidth
\newenvironment{modifyepigraph}[1]
  {\setlength{\epigraphwidth}{#1\textwidth}}
  {\setlength{\epigraphwidth}{\DefaultEpiWidth}}

\usepackage{todonotes}
\begin{document}

\title{Hearts and towers in stable $\infty$-categories}
\author{Domenico Fiorenza, Fosco Loregian, Giovanni Luca Marchetti}
\address{%
Domenico \textsc{Fiorenza} : 
Dipartimento di Matematica ``Guido Castelnuovo'', 
Universit\`a degli Studi di Roma ``la Sapienza'',	
\EmailTo{fiorenza@mat.uniroma1.it}
}
\address{
Fosco \textsc{Loregian} : 
Max Planck Institute for Mathematics,
Bonn, Germany, 
\EmailTo{flore@mpim-bonn.mpg.de}
}
\address{%
Giovanni \textsc{Marchetti} : 
School of Mathematics and Statistics, 
University of Sheffield
\EmailTo{gmarchetti1@shef.ac.uk}
}
\date{\today}
\maketitle

\begin{abstract}
We exploit the equivalence between $t$-structures and normal torsion theories on a stable $\infty$-category to show how a few classical topics in the theory of triangulated categories, i.e., the characterization of bounded $t$-structures in terms of their hearts, their associated cohomology functors, semiorthogonal decompositions, and the theory of tiltings, as well as the more recent notion of Bridgeland's slicings, are all particular instances of a single construction, namely, the tower of a morphism associated with a $J$-slicing of a stable $\infty$-category $\C$, where $J$ is a totally ordered set equipped with a monotone $\mathbb{Z}$-action. 
\end{abstract}

\tableofcontents



\section{Introduction.}\label{abs-intro}
\begin{modifyepigraph}{.5}
\epigraph{If you're going to read this, don't bother.}{C. Palahniuk}
\end{modifyepigraph}

An elementary and yet fundamental theorem in algebraic topology asserts that every sufficiently nice connected topological space $X$ fits into a ``tower''
\[
\cdots \to X_2\to X_1\to X_0=X
\]
 where each $X_n$ is $n$-connected and each map $X_n\to X_{n-1}$ is a fibration that induces isomorphisms in $\pi_{>n}$, and has an Eilenberg-MacLane space $K(\pi_{n}(X),n-1)$ as its fiber.
This result admits an immediate generalization to an arbitrary ambient category which is ``good enough'' for homotopy theory. It is indeed a statement
about the decomposition of an initial morphism $* \to X$ into a tower of fibrations 
whose fibers have homotopy concentrated in a single degree. It is nevertheless only in the setting of $(\infty,1)$-category theory that this result can be given its cleanest conceptualization: the tower of a pointed object $X$ 
is the result of the factorization of $* \to X$ with respect to 
the collection of factorization systems $(n\textsc{-conn}, n\textsc{-trunc})$ whose right classes are given by \emph{$n$-truncated morphisms} \cite[5.2.8.16]{HTT}.

Of course, a similar construction can be exported to \emph{stable} homotopy theory, where the analogue of the factorization system $(n\textsc{-conn}, n\textsc{-trunc})$ is given by the \emph{canonical} $t$-structure $\tee$ on the category of spectra, determined by the objects whose homotopy groups vanish in negative an non-negative degree, respectively,
\begin{gather*}
\Sp_{\ge 0} = \{  A_* \in \Sp \mid \pi_i(A_*)=0;\; i< 0 \}\\
\Sp_{< 0} = \{ B_* \in \Sp \mid \pi_i(B_*)=0;\; i\ge 0 \},
\end{gather*}
together with all its \emph{shifts} $\tee_n=\tee[n]$.\footnote{Here and in the rest of the paper we are implicitly using the equivalence between $t$-structures and \emph{normal torsion theories}: if $\C$ is a stable $\infty$-category with a terminal object, there exists an antitone Galois connection between the poset $\text{Rex}(\C)$ of reflective subcategories of $\C$ and the poset $\pf(\C)$ of prefactorization systems on $\C$ such that $r(\fF)$ is a 3-for-2 class. This adjunction induces a bijective correspondence between the class of certain reflective and coreflective factorization systems called \emph{normal torsion theories} and the class of $t$-structures on (the homotopy category of) $\C$: this statement is the central result of \cite{FL0} where it is called the \emph{Rosetta stone theorem}, and motivates our choice to state our main results in the setting of stable $(\infty,1)$-categories.} In this context, 
it becomes natural to consider the whole $\{\tee_n\mid n\in\mathbb{Z}\}$ as a single object, 
namely the \emph{orbit} of $\tee$ under 
the canonical action of the group of integers on the class $\ts(\Sp)$ of $t$-structures on the category of spectra.
A closer look at this example 
makes it evident that this action is also \emph{monotone} with respect to the natural structure of partially ordered class of $\ts(\Sp)$, and the natural total order of $\mathbb{Z}$: more formally, the group homomorphism $\mathbb{Z}\to \text{Aut}(\ts(\Sp))$ defining the action is also a monotone mapping.

The aim of this paper is to investigate the consequences of taking this point of view further on the classical theory of $t$-structures.  In particular, we describe all the terminology we need about partially ordered groups and their actions in §\refbf{posets}, and then we specialize the discussion to $\Z$\hyp{}actions on partially ordered sets. This is motivated by the fact that the class of $t$-structures on a given stable $\infty$-category carries a natural choice of such an action.

Even if we employ a rather systematic approach, we do not aim at reaching a complete generality, but instead at gathering a number of useful results and nomenclature we can refer to along the present article. Among various possible choices, we mention specialized references as \cite{blyth2005lattices, glass1999partially, Fuch63} for an extended discussion of the theory of actions on ordered groups.

We introduce the definition of a \emph{slicing} of a poset $J$ (\adef\refbf{slicio}) and of a $J$-slicing of a stable $\infty$-category $\C$ (\adef\refbf{J-slicio}) in §\textbf{3}: of course these are not original definitions, as the notion is classical in order theory under different names; our only purpose here is to collect the minimal amount of theory for the sake of clarity. Namely, in the same spirit of Dedekind's construction of real numbers, we consider decompositions of a poset $J$ (more often than not, a totally ordered one) into an \emph{upper} and \emph{lower} set, and associate with each such a decomposition a $t$-structure on an ambient $(\infty,1)$-category $\C$. The totally ordered set $J$ will be assumed to be equipped with a monotone action of $\Z$ and the correspondence
\[
\{\text{slicings of $J$}\} \to \{\text{$t$-structures on $\C$}\}.
\]
will be monotone and $\mathbb{Z}$-equivariant.

Now, following \cite{BBDPervers}, (bounded) $t$-structures on a triangulated category can be seen as the datum of a set of (co)homological functors indexed by integers;  thus a $J$-slicings can be seen as a generalization to ``fractal'' or `` non-integer'' cohomological dimensions now indexed by $J$ and not by $\Z$. Namely, each $J$-slicing induces local cohomology objects depending on an interval, as discussed in §\refbf{sec:towers}. 

In this framework many homological features appear as a shadow of clear constructions with totally ordered sets with $\Z$-actions. For instance, when the totally ordered set $J$ has a heart $J^\heart$, \ie, when there's a $\Z$-equivariant monotone morphism $J \to \mathbb{Z}$, a $J$-slicing on a stable $\infty$-category $\C$ is precisely the datum of a $t$-structure on $\C$ together with a collection of torsion theories parametrized by $J^\heart$ on the heart of $\C$, which turns out to be an abelian $\infty$-category. This is shown §\refbf{hearts} and §\refbf{tiltings}. 

In §\refbf{sec:sods}, we recover the theory of classical semi-orthogonal decompositions by considering the case when $\mathbb{Z}$ acts trivially on $J$. Semi-orthogonal decompositions and $J$-slicing with hearts are essentially the only two interesting classes, as shown by the structure theorem we prove in section §\refbf{concluding}: under suitable finiteness assumptions, the datum of a $J$-slicing on a stable $\infty$-category $\C$ is equivalent to the datum of a finite type semi-orthogonal decomposition of $\C$, together with bounded $t$-structures on the slices and collections of torsion theories  on the hearts of these $t$-structures (Theorem \ref{conclusion}). It is worth mentioning that, under the finiteness conditions of Theorem \ref{conclusion}, when $J=\mathbb{R}$ the notion of $J$-slicing as discussed in the present paper actually becomes a reinterpretation of Bridgeland's definition of slicing of a triangulated category \cite{Brid}. This can be better appreciated by switching to the general approach to `stability data' introduced by \cite{GKR}. 

Finally, in §\refbf{tiltings} we show how the functoriality of the association $J\mapsto \{\text{$J$-slicings}\}$ gives rise to an elegant and synthetic reformulation of classical tilting theory \cite{happel}.
\subsection{Notation and conventions}
We will work within the framework of stable $\infty$-categories in the sense of \cite{LurieHA}. The reader who prefers to work in the more traditional framework of triangulated categories will find no difficulty in pretending that all higher categories are instead categories and that all fiber sequences are distinguished triangles; indeed, many of our statements and constructions are actually adjustments of classical arguments valid in triangulated categories. However a few proofs become more natural when stated in the language of stable $\infty$-categories (for example, theorems whose proof involves a certain universal property unavailable in the triangulated world).

To make the article more self\hyp{}contained and enjoyable to a reader with no previous exposure to stable $\infty$-categories, we recall here the minimal amount of $\infty$-categorical notions we will make use of. As in \cite{HTT}, we will call `$\infty$-category' an $(\infty,1)$-category, i.e. an higher category whose $k$-morphisms are invertible (up to homotopy) for any $k>1$. In colloquial terms, this means that an $\infty$-category has objects, morphisms, homotopies between morphisms, homotopies between homotopies, etc., and such homotopies are all invertible. As shown in \cite{HTT, Joyal2008}, the entire theory of categories transports to $\infty$-categories. In particular, an $\infty$-category $\C$ can have finite co/limits.
\begin{definition*}[Stable $\infty$-category]
An $\infty$-category $\C$ is said to be \emph{stable}  if it has all finite limits and all finite colimits, and if in addition it satisfies the so-called \emph{pullout axiom}: a diagram
\[
\xymatrix{X\ar[r]\ar[d]&Y\ar[d]\\
Z\ar[r]&T
}
\]
in $\C$ (or, more formally, an object of the functor category $\C^\square$) is a pushout if and only if it is a pullback. 
\end{definition*}
We will henceforth call these diagrams \emph{pullout diagrams} or simply \emph{pullouts}.
The pullout axiom is inherently $\infty$-categorical: the only ordinary category with finite limits and colimits satisfying it is the trivial category (i.e. the terminal additive category with a single zero object $\bf 0$). It is also extremely powerful: if $\C$ is a stable $\infty$-category, then its homotopy category $h\C$ is triangulated. In other words, this single and extremely simple axiom subsumes all of the axiomatic of triangulated categories (the notorious octahedral axiom included). 

Unfortunately not every triangulated category can be realized as the homotopy category of a stable $\infty$-category, see \cite{MR2342636}; stable $\infty$-categories therefore only give rise to `well\hyp{}behaved' triangulated categories. It should however be remarked that ill\hyp{}behaved ones are often artificial: the reader can then safely assume that essentially every `reasonable' triangulated category is the homotopy category $h\C$ of some stable $\infty$-category $\C$.

As already said, a particularly pleasant consequence of this good behaviour is the fact that in a stable $\infty$-category the notions of \emph{$t$-structure} \cite{BBDPervers} and of 
\emph{normal factorization system} (or normal torsion theory) \cite{CHK} are naturally equivalent; this remains true in a triangulated category, although the equivalence is much less transparent (see \cite{tderiv}, where this issue is framed in a fairly more general environment). This equivalence will be used several times along the discussion, as well as 
the following result from \cite{FL0}:
\begin{lemma*}[Sator lemma]
If $(\E,\M)$ is a factorization system on a stable $\infty$-category $\C$, with the property that both $\E$ and $\M$ saisfy the `two out of three' property, then for every object $A\in\C$ an initial arrow $0\to A$ lies in $\E$ (resp. in $\M$) \emph{if and only if} the terminal arrow $A\to 0$ of the same object lies in $\E$ (resp. in $\M$).
\end{lemma*}
A final line to conclude this introductory subsection: not to spoil the reader's fun, while avoiding to to hide them their meaning, translations of the quotes opening each section are provided immediately before the bibliography.

\subsubsection*{Acknowledgements} We thank the anonymous referee for useful comments that helped us to improve the overall exposition.

\section{Posets with \texorpdfstring{$\Z$}{Z}-actions.}\label{posets}
\begin{modifyepigraph}{.35}
\epigraph{\japanese{為無為。事無事。味無味。}}{\japanese{老子}}
\end{modifyepigraph}

This section introduces the terminology about partially ordered groups and their actions, and then specializes the discussion to $\Z$-actions. We do not aim at a complete generality, but instead at gathering a number of useful results and nomenclature which is useful to have at hand. Among various possible choices, we mention specialized references as \cite{blyth2005lattices, glass1999partially, Fuch63} for an extended discussion of the theory of actions on ordered groups.

We begin by recalling the well known fact that the category $\Pos$ of partially ordered sets and monotone maps is cartesian closed.
Namely, the product order on the cartesian product $P\times Q$ of two posets $(P,\leq)$ and $(Q,\leq)$, given by $(x,y)\leq (x',y')$ if and only if $x\leq x'$ and $y\leq y'$ makes $(P\times Q,\leq)$ together with the projections on the factors satisfy the universal property of the product in the category of posets. Moreover, the set $\Pos(P,Q)$ of monotone maps from $(P,\leq)$ to $(Q,\leq)$ has a natural order on it given by $f\leq g$ if and only if $f(x)\leq g(x)$ for every $x$ in $P$. With this order one has 
an adjunction
\[
\Pos(P\times Q, R)\cong \Pos(P, \Pos(Q,R)).
\]
\begin{remark}
The product order is not the only standard order one puts on the product $P\times Q$ of two posets $P$ and $Q$. Another commonly used one is the lexicographic order defined by
$(x,y)\leq_{\mathrm{lex}} (x',y')$ if and only if $x< x'$ or $x=x'$ and $y\leq y'$. The lexicographic order does not make $P\times Q$ be the product of $P$ and $Q$ in the category of posets, but it still has a few peculiar properties that, as we are going to see, are relevant to the theory of slicings.
\end{remark}

\begin{remark}
\label{quotient1}
Let $(P,\leq)$ be a poset, and let $\sim$ be an equivalence relation on the set $P$. One says that $\sim$ is compatible with the order relation if $x\leq y$, $x\sim x'$ and $y\sim y'$ imply $x'\leq y'$ or $x'\sim y'$. When this happens the quotient set $P/_{\!\sim}$ inherits a order relation from $P$ by $[x]\leq [y]$ if and only if $x\leq y$ or $x\sim y$. Moreover the projection to the quotient
$P\to P/_{\!\sim}$ is a morphism of posets.
\end{remark}

\begin{definition}
A \emph{partially ordered group} (``po-group'' for short) is a pair $(G, \leq)$ consisting of a group $G$ and of a partial order relation $\leq $ on $G$ such that the group multiplication  $\cdot \colon G\times G\to G$ is a map of posets, where $G\times G$ is endowed with the product order: for any two pairs $(g,h)$ and $(g',h')$ with $g\leq g'$ and $h\leq h'$ we have $gg' \le hh'$.
\end{definition}
\begin{remark}
In the literature on the subject it is customary draw a distinction between a \emph{left} po-group and a \emph{right} po-group. We choose to ignore this subtlety, since all the po-groups we will be dealing with will be ordered by two-sided congruences.
\end{remark}
\begin{remark} If $(G,\leq)$ is a po-group, the inversion $(\firstblank)^{-1}\colon G \to G$ is an antitone antiautomorphism of groups, \ie we have that
$g\leq h$ if and only if $h^{-1}\leq g^{-1}$. Moreover
the set $G^+$ of \emph{positive} elements, i.e. the set $\{g\in G\mid 1\leq  g\}$ is closed under conjugation.
\end{remark}
\begin{example}
Let $(P,\leq)$ be a poset, and let $\mathrm{Aut}_{\Pos}(P)$ be the automorphism group of $P$ as a poset, i.e., the set of monotone bijections of $p$ into itself. Then $\mathrm{Aut}_{\Pos}(P)$ inherits an order relation by its inclusion in the poset $\Pos(P,P)$, and this makes $\mathrm{Aut}_{\Pos}(P)$ a partially ordered group. This is the standard po-group structure on $\mathrm{Aut}_{\Pos}(P)$.
\end{example}
\begin{remark}
Any group $G$ can be seen as a po-group with the trivial order relation $g\leq h$ if and only if $g=h$. It is worth noticing that on finite groups the trivial order is the only possible po-group structure. Namely, assume $g\leq h$ and let $k=g^{-1}h$. Then we have $1\leq k\leq k^2\leq k^2\leq\dots\leq k^{\mathrm{ord}(k)}=1$ and so $k=1$, i.e., $g=h$.
\end{remark}
\begin{definition}
A \emph{homomorphism} of po-groups consists of a group morphism $f\colon G\to H$ which is also a monotone mapping. This 
defines a category $\cate{PoGrp}$ of partially ordered groups and their homomorphisms.
\end{definition}
\begin{definition}\label{def:g.poset}
Let $(G,\leq)$ be a po-group. A $G$-poset is a partially ordered set $(P, \le)$ endowed with a po-group homomorphism $G\to \text{Aut}_{\Pos} (P)$ to the group of order isomorphisms of $P$ with its standard po-group structure.
\end{definition}
\begin{remark}
Equivalently, a $G$-poset is a partially ordered set $P$ together with a group action $G\times P\to P$ which is a morphism of posets, where on $G\times P$ one has the product order.
\end{remark}
\begin{example}
Every po-group $G$ is a $G$-poset with the multiplication action of $G$ on itself. 
\end{example}
\begin{remark}\label{quotient2}
An equivalence relation $\sim$ on a $G$-poset $P$ is said to be compatible with the $G$-action if $x\sim y$ implies $g\cdot x\sim g\cdot y$ for any $g$ in $G$. If $\sim$ is compatible both with the order and with the $G$-action then the quotient set $P/_{\!\sim}$ is naturally a $G$-poset with the $G$-action $g\cdot[x]=[g\cdot x]$. Moreover the projection to the quotient is a morphism of $G$-posets.
\end{remark}

We now specialize our discussion to the case $G=\Z$
\begin{definition}\label{zposet}
A $\Z$-\emph{poset} is a partially ordered set $(P,\leq)$ together with a group action 
\[
+_P\colon P\times \Z \to P \colon (x,n) \mapsto x+_Pn 
\]
 which is a morphism of partially ordered sets, when $\Z$ is regarded with its usual total order.
\end{definition}
\begin{remark}\label{trivial.but.useful}
It is immediate to see that a $\Z$-poset is equivalently the datum of a poset $(P,\leq)$ together with a monotone bijection $\rho\colon P\to P$ such that $x\leq \rho(x)$ for any $x$ in $P$. The function $\rho$ and the action are related by the identity $\rho(x)=x+_P1$.
\end{remark}
\begin{notat}
To avoid a cumbersome accumulation of indices, the action $+_P$ will be often denoted as a simple ``$+$''. This is meant to evoke in the reader the most natural example of a $\Z$-poset, given by $\mathbb{Z}$ itself, and to simplify our notation for the axioms of an action:
\[
\begin{cases}
(x +_P m) +_P n = x +_P (m+n);\\
x +_P 0 = x.
\end{cases}
\]
We will also write $x-n$ for $x+_P(-n)$.
\end{notat}
\begin{example}The poset $(\Z ,\leq)$ of integers with their usual order is a $\Z$-poset with the action given by the usual sum of integers. The poset $(\R,\leq)$ of real numbers with their usual order is a $\Z$-poset for the action given by the sum of real numbers with integers (seen as a subring of real numbers).
\end{example}
\begin{example}
Given any poset $(P,\leq)$, the poset $\Z\times_{\mathrm{lex}}P$ carries a natural $\mathbb{Z}$-action given by $(n,x)+1=(n+1,x)$, i.e., by the standard $\Z$-action on the first factor and by the trivial $\Z$-action on the second factor. 
\end{example}

\begin{remark}\label{rem.finite}
 If $(P,\leq)$ is a finite poset, then the only $\Z$-action it carries is the trivial one. Indeed, if $\rho\colon P\to P$ is the monotone bijection associated with the $\Z$-action, one sees that $\rho$ is a finite order element in $\mathrm{Aut}_{\Pos}(P)$, by the finiteness of $P$. Therefore there exists an $n\geq 1$ such that $\rho^n=\mathrm{id}_P$. It follows that, for any $x$ in $P$,
 \[
 x\leq x+1\leq\cdots\leq x+n=x
 \]
and so $x=x+1$.
\end{remark}
\begin{remark}\label{minmax}
An obvious terminology: a \emph{$G$-fixed point} for a $G$-poset $P$ is an element $p\in P$ kept fixed by all the elements of $G$ under the $G$-action. 
Clearly, an element $x$ of a $\Z$-poset $P$ is a $\Z$-fixed point if and only if $x+ 1 = x$, or equivalently $x-1=x$. From this it immediately follows that if $x\in P$ is a $\le$-maximal or $\le$-minimal element in the $\Z$-poset $P$, then it is a $\Z$-fixed point.
\end{remark}
\begin{remark}
Given a poset $P$ we can always define a partial order on the set $P_{\!\bowtie} = P\cup\{-\infty,+\infty\}$ which extends the partial order on $P$ by the rule $-\infty\leq x\leq +\infty$ for any $x\in P$. 
\end{remark}
\begin{lemma}
 If $(P,\leq)$ is a $\Z$-poset, then $(P_{\!\bowtie},\leq)$ carries a natural $\Z$-action extending the $\Z$-action on $P$, by declaring both $-\infty$ and $+\infty$ to be $\Z$-fixed points.
\end{lemma}
\begin{proof}
 Adding a fixed point always gives an extension of an action, so we only need to check that the extended action is compatible with the partial order. This is equivalent to checking that also on $P_{\!\bowtie}$ the map $x\mapsto x+1$ is a monotone bijection such that $x\leq x+1$, which is immediate. 
\end{proof}
Posets with $\Z$-actions naturally form a category $\ZPos$, whose morphisms are \emph{$\Z$-equivariant} morphisms of posets. More explicitly, if $P$ and $Q$ are $\Z$-posets with actions $+_P$ and $+_Q$, then a morphism of $\Z$-posets between them is a morphism of posets $\varphi\colon P\to Q$ such that
\[
\varphi(x+_P n)=\varphi(x)+_Q n,
\]
for any $x\in P$ and any $n\in \Z$.
\begin{remark}
If $P$ and $Q$ are $\Z$-posets, then the hom-set $\ZPos(P,Q)$ is naturally a $\Z$-poset. Namely, as we have already remarked, $\Pos(P,Q)$ is naturally a poset, and so $\ZPos(P,Q)$ inherits the poset structure as a subset. The $\Z$-action is given by $(\varphi+n)(x)=\varphi(x) +_Q n$. This makes $\ZPos$ a closed category.
\end{remark}

\begin{remark}
Every poset can be seen as a $\Z$-poset with the trivial $\Z$-action. Since every monotone mapping is $\Z$-equivariant with respect to the trivial $\Z$-action, this gives a fully faithful embedding $\Pos\to \ZPos$.
\end{remark}
\begin{lemma}\label{trivial.but.useful2}
The choice of an element $x$ in a $\Z$-poset $P$ is equivalent to the datum of a $\Z$-equivariant morphism $\varphi\colon(\Z ,\leq)\to (P,\leq)$. Moreover $x$ is a $\Z$-fixed point if and only if the corresponding morphism $\varphi$ factors $\Z$-equivariantly through $(*,\leq)$, where $*$ denotes the terminal object of $\Pos$. 
 \end{lemma}
\begin{proof}
To the element $x$ one associates the $\Z$-equivariant morphism $\varphi_x$ defined by $\varphi_x(n)=x+n$. To the $\Z$-equivariant morphism $\varphi$ one associates the element $x_\varphi=\varphi(0)$. It is immediate to check that the two constructions are mutually inverse. The proof of the second part of the statement is straightforward.
\end{proof}
\begin{lemma}
Let $\varphi\colon(\Z ,\leq)\to (P,\leq)$ be a $\Z$-equivariant morphism of $\Z$-posets. Then $\varphi$ is either injective or constant.
\end{lemma}
\begin{proof}
Assume $\varphi$ is not injective. then there exist two integers $n$ and $m$ with $n>m$ such that $\varphi(n)=\varphi(m)$. By $\Z$-equivariancy we therefore have
\[
x_\varphi+(n-m)=x_\varphi,
\]
with $n-m\geq 1$ and $x_\varphi=\varphi(0)$. The conclusion then follows by the same argument used in Remark \refbf{rem.finite}.
\end{proof}
\begin{lemma}\label{extends}
Let $\varphi\colon (P,\leq)\to (Q,\leq)$ be a morphism of $\Z$-posets. Assume $Q$ has a minimum and a maximum. Then $\varphi$ extends to a morphism of $\Z$-posets $(P_{\!\bowtie},\leq)\to (Q,\leq)$ by setting $\varphi(-\infty)=\min(Q)$ and $\varphi(+\infty)=\max(Q)$.
\end{lemma}
\begin{proof}
Since $\min(Q)$ and $\max(Q)$ are $\Z$-fixed points by Remark \refbf{minmax}, the extended $\varphi$ is a morphism of $\Z$-posets. Moreover, since $\min(Q)$ and $\max(Q)$ are the minimum and the maximum of $Q$, respectively, the extended $\varphi$ is indeed a morphism of posets, and so it is a morphism of $\Z$-posets.
\end{proof}

All of the above applies in particular to totally ordered sets. We will denote by $\Tos\subseteq \Pos$ the full subcategory of totally ordered sets, and by $\mathbb{Z}\text{-}\Tos\subseteq \ZPos$ the full subcategory of $\Z$-actions on totally ordered sets.

\begin{lemma}\label{equivalence}
Let $(P, \leq)$ be a totally ordered $\mathbb{Z}$-poset. The relation $x\sim y$ if and only if there are integers $a,b \in \mathbb{Z}$ such that $x + a \le y \le x + b $ is an equivalence relation on $P$ compatible with both the order and the $\mathbb{Z}$-action. It therefore induces a morphism of $\mathbb{Z}$-posets $P\to P/_{\!\sim}$ given by the projection to the quotient. Moreover, $P/_{\!\sim}$ is totally ordered and the $\mathbb{Z}$-action on the quotient is trivial.
\end{lemma}
\begin{proof}
Checking that $\sim$ is an equivalence relation is immediate: reflexivity is manifest; symmetry reduces to noticing that $x + a \le y \le x + b $ is equivalent to $y-b\leq x\leq y-a$; transitivity follows by  the fact that $x + a \le y \le x + b $ and $y + c \le z \le y + d$ together imply $x+ (a+ c) \le z \le x+(b + d)$. To see that $\sim$ is compatible with the order relation, let $x\leq y$ and let $x\sim x'$ and $y\sim y'$. Then there exist $a, b, c$ and $d$ in $\mathbb{Z}$ such that $x + a \le x' \le x + b $ and $y + c \le y' \le y + d$. Since $P$ is totally ordered, either $x'\leq y'$ or $y'\leq x'$. In the second case we have $y'\leq x'  \le x + b\le y+b \le y'+b-c$, and so $x'\sim y'$ by definition of the relation $\sim$. The compatibility of $\sim$ with the $\mathbb{Z}$-action is straightforward. Therefore by Remarks \refbf{quotient1} and \refbf{quotient2} we see that the projection to the quotient $P\to P/_{\!\sim}$ is a morphism of $\mathbb{Z}$-posets. Since the order on $P$ is total, so is also the order induced by $\sim$ on the quotient set. Finally, to see that the $\mathbb{Z}$ action on $P/_{\!\sim}$ is trivial, just notice that for any $x$ in $P$, we have $x\leq x+1\leq x+1$ and so $[x]+1=[x]$.
\end{proof}

\begin{remark}\label{iq=q}
If the $\mathbb{Z}$-action on the totally ordered set $P$ is trivial, then the equivalence relation from Lemma \refbf{equivalence} is trivial as well: $x\sim y$ if and only if $x=y$.
\end{remark}
\begin{lemma}\label{lemma.representatives}
Let $(P, \leq)$ be a totally ordered $\mathbb{Z}$-poset, and let $\sim$ be the equivalence relation from Lemma \refbf{equivalence}. Then either $[x]=\{x\}$ or $[x]=\mathbb{Z}\times_{\mathrm{lex}}[x,x+1)$, where $[x,x+1)=\{y\in P\mid x\leq y<x+1\}$
\end{lemma}
\begin{proof}
Let $x\in P$; then either $x=x+1$ or $x<x+1$. In the first case $x$ is a fixed point for the $\Z$-action on $P$ and so the equivalence relation $\sim$ is the trivial one: $y\sim x$ if and only if $y=x$. If $x<x+1$ then the interval $[x,x+1)$ is nonempty, as $x\in [x,x+1)$.

Let $\varphi\colon \mathbb{Z}\times_{\mathrm{lex}}[x,x+1)\to P$ the map defined by $(n,y)\mapsto y+n$. The map $\varphi$ is a morphism of $\Z$-posets. 

Indeed, if $(n,y)\leq_{\mathrm{lex}} (n',y')$ either $n<n'$ or $n=n'$ and $y\leq y'$. In the first case we have $n+1\leq n'$ and so $y+n< x+1+n\leq x+n'\leq y'+n'$, whereas in the second case we have $y+n\leq y'+n=y'+n'$. The map $\varphi$ is also injective. 

To conclude we only need to show that $\varphi\colon  \mathbb{Z}\times_{\mathrm{lex}}[x,x+1)\to [x]$ is surjective. Let $z\in [x]$. By definition of the equivalence relation there exist $a,b$ in $\mathbb{Z}$ such that $x+a\leq z\leq x+b$. Since $x$ is not a $\Z$-fixed point, we have $a\leq b$ and $z\in [x+a,x+b+1)$. Writing
\[
[x+a,x+b+1)=\bigcup_{n=a}^{b}[x+n,x+n+1)
\]
we see that there exists $n\in \Z$ such that $z\in [x,x+1)+n$, i.e., $z=\varphi(n,y)$ for some $y$ in $[x,x+1)$.
\end{proof}
\begin{proposition}\label{adjoint}
The fully faithful embedding $(\firstblank)^\flat\colon \Tos\to \Z\text{-}\Tos$ given by trivial $\Z$-actions on totally ordered sets has a left adjoint.
\end{proposition}
\begin{proof}
For any totally ordered $\Z$-poset $P$, let $\iota(P)$ be the $\mathbb{Z}$-poset $P/_{\!\sim}$, where $\sim$ is the equivalence relation from Lemma \refbf{equivalence}. Then $I\colon \Z\text{-}\Tos\to \Tos$ is a functor, since if $f\colon P\to Q$ is a morphism of $\mathbb{Z}$-posets then $x+a\leq y\leq x+b$ implies
\[
f(x)+a=f(x+a)\leq f(y)\leq f(x+b)=f(x)+b
\]
and so $f$ induces a well defined morphism of sets $\tilde{f}\colon \iota(P)\to \iota(Q)$. It is immediate to see that $\tilde{f}$ is actually a morphism of $\Z$-posets and that $f\rightsquigarrow \tilde{f}$ preserves identities and compositions of morphisms. Finally, to see that $I$ is a right adjoint to the trivial action embedding $\Tos\to \Z\text{-}\Tos$, let $P$ be a a totally ordered $\mathbb{Z}$-poset and $Q$ be a totally ordered set. Since $I$ is a functor, a morphism of $\mathbb{Z}$-posets $f\colon P\to Q^\flat$ induces a morphism of posets $\tilde{f}\colon \iota(P)\to \iota(Q^\flat)$. Moreover, since the $\mathbb{Z}$-action on $Q^\flat$ is trivial, we have $\iota(Q^\flat)\cong Q$, see Remark \refbf{iq=q}. Therefore, $f\mapsto \tilde{f}$ defines a map
\[
\Z\text{-}\Tos(P,Q^\flat)\to \Tos(\iota(P),Q)
\]
which we want to show is a bijection. Assume $\tilde{f}_1=\tilde{f}_2$. Then, for any $x$ in $P$ we have $f_1(x)\sim f_2(x)$ in $Q^\flat$. Since the equivalence relation on $Q^\flat$ is trivial, this means $f_1=f_2$, so $f\rightsquigarrow \tilde{f}$ is injective. Let now $\varphi$ be a morphism of posets, $\varphi\colon \iota(P)\to Q$, and let $f=\varphi\circ \pi$ where $\pi\colon P\to \iota(P)$ is the projection to the quotient. We have $\tilde{f}([x])=[f(x)]=f(x)=\varphi(\pi(x))=\varphi([x])$, and so $\varphi=\tilde{f}$, i.e., $f\rightsquigarrow \tilde{f}$ is surjective.
\end{proof}

\section{Histoire d'\texorpdfstring{$\O(J)$}{\O(J)}.}\label{histoire}
\begin{modifyepigraph}{.57}
\epigraph{Sa liberté était pire que n'importe quelle chaîne}{D\@. Aury}
\end{modifyepigraph}
Recall that a lower set in a poset $J$ is a subset $L\subseteq J$ such that if $x\in L$ and $y\leq x$ then $y\in L$; the set of lower sets of $J$ is denoted ${\downharpoonleft} J$ and it naturally a partially ordered set. Dually, one defines upper sets and the set ${\upharpoonleft}J$ of upper sets with its natural partial order.
\begin{definition}\label{slicio}
Let $J$ be a poset. A \emph{slicing} of $J$ is a pair $(L,U)$, where $L$ is a lower set in $J$, $U$ is an upper set, $L\cap U=\emptyset$ and $L\cup U=J$. The collection of all slicings of $J$ will be denoted by $\O(J)$.
\end{definition}

\begin{remark}
Since the complement of an upper set is a lower set and vice versa, the projection on the second factor is a bijection
\[
\O(J)\xrightarrow{\sim} {\upharpoonleft}J
\]
This induces a natural partial order on $\O(J)$: we set $(L_1,U_1)\leq (L_2,U_2)$ if and only if $U_2\subseteq U_1$. Notice that $\O(J)$ has a minimum given by the slicing $(\emptyset,J)$ and a maximum given by the slicing $(J,\emptyset)$.
\end{remark}

\begin{remark}\label{rem.totOj}
If $J$ is totally ordered, then so is $\O(J)$. Namely, let $U_1$ and $U_2$ two upper sets in $J$ and assume that $U_1$ is not a subset of $U_2$. Then there exists an element $x$ in $U_1$ which is not in $U_2$. If $y\in U_2$ then either $y\leq x$ or $y\geq x$ since $J$ is totally ordered. But since $U_2$ is an upper set $y\leq x$ would imply $x\in U_2$ against our assumption. This means that $y\geq x$ and, since $U_1$ is an upper set, this implies $y\in U_1$. Therefore if $U_1$ is not a subset of $U_2$ we have that $U_2\subseteq U_1$.
\end{remark}
\begin{remark}
If $J$ is a $\mathbb{Z}$-poset, then so is $\O(J)$. The natural $\mathbb{Z}$-action on $\O(J)$ is given by
\[
(L,U)+n=(L+n,U+n),
\]
where $L+n=\{x+n\,|\, x\in L\}$ and $U+n=\{x+n\,|\, x\in U\}$.
\end{remark}
\begin{remark}
Every element $x$ in $J$ determines two slicings of $J$: $((-\infty,x),[x,+\infty))$ and $((-\infty,x],(x,+\infty))$. Here $(-\infty,x)$ is the lower set $\{y\in J\,|\, y<x\}$, and similarly for $(-\infty,x]$, $(x,+\infty)$ and $[x,+\infty)$. This gives two natural morphisms of posets $J\to \O(J)$. If $J$ is a $\mathbb{Z}$-poset, then these morphisms are $\mathbb{Z}$-equivariant.
\end{remark}
The construction of $\O(J)$ is actually functorial in $J$ so that we have the following
\begin{lemma}\label{lemma.O-is-functor}
The map $J\rightsquigarrow \O(J)$ defines a functor
\[
\O\colon \ZPos^{\mathrm{op}} \to \ZPos_{\top}
\]
where $\ZPos_{\top}$ denotes the category of $\Z$-posets with minimum and maximum and with $\Z$-morphism of posets preserving them (these maps are called \emph{$\{0,1\}$-homomorphisms} in lattice theory, see \cite{Gratzer}).
\end{lemma}
\begin{proof}
By the above remarks the only thing we have to prove is functoriality. For any morphism of $\Z$-posets $f\colon J_1\to J_2$, we set $\O(f)\colon U\to f^{-1}(U)$. It is immediate to see that $f^{-1}(U)$ is an upper set in $J_2$ for any upper set $U$ in $J_1$ an that $f^{-1}(U+1)=f^{-1}(U)+1$, so that $\O(f)$ is indeed a morphism of $\Z$-posets from $\O(J_2)$ to $\O(J_1)$. Moreover we have $\O(\mathrm{id}_J)=\mathrm{id}_{\O(J)}$ and $\O(fg)=\O(g)\O(f)$, and $\O(f)(\emptyset)=\emptyset$ and $\O(f)(J_2)=J_1$.
\end{proof}
\begin{remark}\label{rem.hom-is-Z-pos}
Since the minimum and the maximum of a $\Z$-poset, when they exist, are necessarily fixed points of the $\Z$-action, we see that the inclusion  $\ZPos_{\top}(P,Q)\subseteq \ZPos(P,Q)$ induces a $\Z$-poset structure on $\ZPos_{\top}(P,Q)$ making this inclusion a morphism of $\Z$-posets.
\end{remark}

\begin{example}\label{ex.Z-and-R}
The morphism $n\mapsto [n,+\infty)$ induces an isomorphism of $\mathbb{Z}$-posets $\mathbb{Z}_{\bowtie}\xrightarrow{\sim} \O(\mathbb{Z})$. The morphisms $x\mapsto (x,+\infty)$ and $x\mapsto [x,+\infty)$ together induce an isomorphism of posets $\mathbb{Z}$-posets
\[
(\mathbb{R}\times_{\mathrm{lex}} \uno)_{\bowtie}\xrightarrow{\sim} \O(\mathbb{R}),
\]
where $\uno $ is the totally ordered set $\{0 < 1\}$.
\end{example}
\begin{definition}
Before introducing the main definition of this section, let us recall that a \emph{$t$-structure} on a stable $\infty$-category $\C$ consists of a pair $\tee=(\cate{L},\cate{U})$ of full sub-$\infty$-categories satisfying the following properties:
\begin{enumerate}[label=$\roman*$)]
\item orthogonality: $\C(X, Y)$ is  contractible for each $X\in \cate{U}$, $Y\in \cate{L}$;
\item one has $\cate{U}[1]\subseteq \cate{U}$ and $\cate{L}[-1]\subseteq \cate{L}$;
\item Any object $X\in\C$ fits into a (homotopy) fiber sequence $X_{\cate{U}}\to X\to X_{\cate{L}}$, with $X_{\cate{U}}$ in $\cate{U}$ and $X_{\cate{L}}$ in $\cate{L}$. 
\end{enumerate}
\end{definition}
We introduce further terminology as a separate remark:
\begin{remark}
The categories $\cate{L}$ and $\cate{U}$ are called the \emph{lower sub-$\infty$-category} and the \emph{upper sub-$\infty$-category} of the $t$-structure $\tee$, respectively. 

The collection $\ts(\C)$ of all $t$-structures on a stable $\infty$-category $\C$ is a poset with respect to following order relation: given two $t$-structures $\tee_1=(\cate{L}_1, \cate{U}_1)$ and  $\tee_2=(\cate{L}_2, \cate{U}_2)$, one has  $\tee_1 \leq \tee_2$ if and only if $\cate{U}_2\subseteq \cate{U}_1$. 

The ordered group $\Z $ acts on $\ts(\C)$ in a way that is fixed by the action of the generator $+1$; this maps a $t$-structure $\tee=(\cate{L},\cate{U})$ to the \emph{shifted} $t$-structure $\tee[1]=(\cate{L}[1],\cate{U}[1])$. Since $\tee\leq\tee[1]$ one sees that $\ts(\C)$ is naturally a $\Z $-poset. Finally, the poset $\ts(\C)$ has a minimum and a maximum given by $(\mathbf{0},\C)$ and $(\C,\mathbf{0})$, respectively. These are called the \emph{trivial} $t$-structures.
\end{remark}
\begin{definition}\label{J-slicio}
Let $(J,\leq)$ be a $\Z $-poset. A \emph{$J$-slicing} of a stable $\infty$-category $\C$ is a $\Z $-equivariant morphism of posets $\tee\colon \O(J)\to \ts(\C)$ respecting minima and maxima on both sides. We denote as $\slicings(J,\C)$ the class of all $J$-slicings of the category $\C$;\footnote{The Japanese verb {\japanese[gbsn]{切る}} (``kiru'', \emph{to cut}) contains the radical {\japanese{切}}, the same of \emph{katana}.}
 \end{definition}
 More explicitly, a $J$-slicing is a family $\{\tee_\lu\}_{(L,U)\in \O(J)}$ of $t$-structures on $\C$ such that
 \begin{enumerate}[label=$\roman*$)]
\item $\tee_{(L_1,U_1)}\leq \tee_{(L_2,U_2)}$ if $(L_1,U_1)\leq (L_2,U_2)$ in $\O(J)$;
\item $\tee_{(L,U)+1}=\tee_\lu[1]$ for any $(L,U)\in \O(J)$.
\item $\tee_{(\emptyset,J)}=(\mathbf{0},\C)$ and $\tee_{(J,\emptyset)}=(\C,\mathbf{0})$.
\end{enumerate}
\begin{remark}\label{rem.slicing-functor}
Of course, $\slicings(J,\C)=\ZPos_{\top}(\O(J),\ts(\C))$: this, together with \refbf{lemma.O-is-functor} and \refbf{rem.hom-is-Z-pos}, gives that $J\mapsto \slicings(J,\C)$ is a functor.
\end{remark}

\begin{notat}\label{magictrick}
We will denote the lower and the upper sub-$\infty$-categories of the $t$-structure $\tee_\lu$ by $\C_L$ and $\C_U$, respectively, i.e., we write $\tee_\lu=(\C_L,\C_U)$.
For $i\in J$, we will write $\C_{\geq i}$,  $\C_{> i}$, $\C_{\leq i}$  and $\C_{<i}$ for $\C_{[i,+\infty)}$, $\C_{(i,+\infty)}$, $\C_{(-\infty,i]}$ and $\C_{(-\infty,i)}$, respectively.  Note that, by $\Z $-equivariancy, we have $
\C_{\geq i+1}=\C_{\geq i}[1]$, and similarly for the other cases.
\end{notat}

\begin{example}\label{ex.Z-is-t}
By Lemma \refbf{trivial.but.useful2} and Example \refbf{ex.Z-and-R}, a $\Z $-slicing on $\C$ is equivalent to the datum of a $t$-structure $\tee_0=(\C_{<0},\C_{\geq 0})$. One has $\tee_n=(\C_{<n},\C_{\geq n})$ for any $n\in \mathbb{Z}$, consistently with the Notation \refbf{magictrick}, $\tee_{-\infty}=(\mathbf{0},\C)$ and $\tee_{+\infty}=(\C,\mathbf{0})$. Notice that by our Remark \refbf{rem.finite}, as soon as $\C_{\ge 1}$ is a proper subcategory of $\C_{\ge 0}$, then the inclusion $\C_{\geq n+1}\subseteq \C_{\geq n}$ is proper for all $n\in\Z$, \ie the orbit $\tee + \Z$ is an infinite set. The equivalence between $t$-structures and $\Z$-slicings can also be seen in the light of Remark \refbf{rem.slicing-functor}: for every $\Z$-poset $P$ with minimum and maximum one has a distinguished isomorphism $\ZPos_{\top}(\mathcal{O}(\Z),P)\xto{\sim}P$ given by $\varphi\mapsto \varphi([0,+\infty))$
\end{example}
\begin{example}\label{what.s.slici}
By Example \refbf{ex.Z-and-R}, an $\R$-slicing on $\C$ is the datum of two $t$-structures $(\C_{<\lambda},\C_{\geq \lambda})$ and  $(\C_{\leq \lambda},\C_{> \lambda})$ on $\C$ for any $\lambda\in \R$ in such a way that $\C_{\geq \lambda+1}=\C_{\geq \lambda}[1]$, etc., and with the inclusions $\C_{>\lambda}\subseteq \C_{\geq\lambda}$ for any $\lambda\in \mathbb{R}$ and 
\[
\C_{>\lambda_2}\subseteq \C_{\geq \lambda_2} \subseteq \C_{>\lambda_1}\subseteq \C_{\geq \lambda_1}
\]
for any $\lambda_1<\lambda_2$ in $\mathbb{R}$. $\mathbb{R}$-slicings have been introduced in \cite{Brid}, where they are called simply ``slicings''. Actually \cite{Brid} imposes more restrictive conditions to ensure ``compactness'' of the factorization, we will come back to this later. Compare also \cite{GKR}.
\end{example}

\begin{remark}
Since the subcategories $\C_L$ and  $\C_U$ are the lower and the upper subcategories of a $t$-structure $\tee_\lu$ they are reflexive and coreflective, respectively. In particular we have reflection and coreflection functors
\[
R_L\colon \C\to \C_L;\qquad\qquad S_U\colon \C\to \C_U.
\]
For $X$ an object in $\C$ we will occasionally write $X_L$ for $R_LX$ and $X_U$ for $S_UX$, and similarly for morphisms. Finally, by composing $R_L$and $S_U$ with the inclusions of $\C_L$ and $\C_U$ in $\C$, we can look at $R_L$ and $S_U$ as endofunctors of $\C$. 
\end{remark}

In order to investigate properties of the (co-)reflections $S$ and $R$, we recall the main result from \cite{FL0}: there is an equivalence between $t$-structures on $\C$ and normal factorization systems on $\C$, so that $\tee$ can equivalently be seen as a $\Z$-equivariant morphism $\O(J)\to  \fs(\C)$, where $\fs(\C)$ denotes the $\mathbb{Z}$-poset of normal factorization systems of $\C$. Explicitly, this equivalence is given as follows: given a $t$-structure $(\cate{L},\cate{U})$ on $\C$,
the corresponding factorization system $(\E,\M)$ is characterized by
\begin{align*}
0\to X &\text{ is in }\E\text{ if and only if }X\in \cate{U}\\
X\to 0 &\text{ is in }\M\text{ if and only if }X\in \cate{L}
\end{align*}
Since we are going to use this fact several times, we recall that both the class $\E$ and the class $\M$ have the 3-for-2 property. In particular this implies the Sator lemma:
\begin{align*}
(0\to X) &\text{ is in }\E \text{ if and only if }(X\to 0) \text{ is in }\E\\
(X\to 0) &\text{ is in }\M \text{ if and only if }(0\to X) \text{ is in }\M
\end{align*}
For further information on normal factorization systems in stable $\infty$-categories we address the reader to \cite{FL0,tstructures}.

\begin{remark}
Notice that the left class $\E$ of the normal factorization system $(\E,\M)$ corresponds to the right class $\cate{U}$ of the corresponding $t$-structure $(\cate{L},\cate{U})$. One could avoid this position switch by writing the pair of classes in a $t$-structure as $(\cate{U},\cate{L})$, however we preferred to keep the upper class on the right to agree with the standard orientation on the line of real numbers.
\end{remark}

\begin{remark}\label{closed.by.extensions} Since $(\E,\M)$ are a factorization system, the class $\M$ is closed under pullbacks and the class $\E$ is closed under pushouts. Together with the Sator lemma this implies that  $\cate{L}$ and $\cate{U}$ are extension closed.

\begin{lemma}
\label{lem.iterated} 
Let $\tee$ be a $J$-slicing of $\C$ and let $(L_0,U_0)$ and $(L_1,U_1)$ two slicings of $J$ with $(L_0,U_0)\leq (L_1,U_1)$. Then we have the natural isomorphisms
\begin{enumerate}[label=$\roman*$)]
\item $R_{L_0}R_{L_1}\cong R_{L_1}R_{L_0}\cong R_{L_0}$;
\item $ S_{U_0}S_{U_1}\cong S_{U_1}S_{U_0}\cong S_{U_1}$;
\item $R_{L_0}S_{U_1}\cong S_{U_1}R_{L_0}\cong 0$.
\end{enumerate}
\end{lemma}
\begin{proof} 
It is enough to prove $(i)$ and $(iii)$, as the proof of $(ii)$ is dual to $(i)$.

We denote $\tee_0$ and $\tee_1$ the $t$-structures corresponding to $(L_0,U_0)$ and $(L_1,U_1)$, respectively, and by $(\E_0,\M_0)$ and $(\E_1,\M_1)$ the corresponding normal torsion theories. Since $(L_0,U_0)\leq (L_1,U_1)$ we have $\M_{0} \subseteq \M_{1}$ and $\E_1\subseteq \E_0$.

The object $S_1X$ is obtained by $(\E_1, \M_1)$-factoring the arrow $0\to X$ as $0\xto{\E_1}S_1X\xto{\M_1}X$. Since $\E_1\subseteq\E_0$, this shows that $0\to S_1X$ is in $\E_0$ and so, by the Sator lemma, also $S_1X\to 0$ is in $\E_0$. Now the object $R_0S_1X$ is obtained by the $(\E_0, \M_0)$-factorization of the terminal morphism $S_1X\to 0$ as $S_1X\xto{\E_0}R_0S_1X\xto{\M_0}0$.

By the 3-for-2 property for $\E_0$ we see that $R_0S_1X \to 0$ lies in $\E_{0}\cap \M_{0}$, hence it an isomorphism, so that $R_0S_1 X\cong 0$. The proof that $S_1R_0X\cong X$ is perfectly dual. This proves $(iii)$.

The proof of $(i)$ goes as follows. The reflection $R_0R_1 X$ is defined by the $(\E_0,\M_0)$-factorization $R_1X\xto{\E_0}R_0R_1X\xto{\M_0}0$. Since $\E_1\subseteq \E_0$ and $X\to R_1X$ is in $\E_1$ we see that $X\to R_0R_1X\to 0$ is already the $(\E_0,\M_0)$-factorization of $X\to 0$ and so by uniqueness of the factorization we have $R_0R_1X\cong R_0X$. Finally, the reflection $R_1R_0 X$ is defined by the $(\E_1,\M_1)$-factorization $R_0X\xto{\E_1}R_1R_0X\xto{\M_1}0$. Since $R_0X\to 0$ is in $\M_0\subseteq \M_1$, by the 3-for-2 property we have that $R_0X\to R_1R_0X$ is in $\E_1\cap \M_1$ and so it is an isomorphism. 
\end{proof}
\begin{remark}
Notice how, in the proof of the above lemma, one sees that applying $R_{L_0}$ to the natural morphism $R_{L_1}X\to 0$ we get a natural morphism $R_{L_1}X\to R_{L_0}  X$, an so one has a natural transformation $R_{L_1}\to R_{L_0}$. Dually, we have a natural transformation $S_{U_1}\to S_{U_0}$.
\end{remark}
\begin{lemma}\label{lem.for-homology} 
Let $\tee$ be a $J$-slicing of $\C$ and let $(L_0,U_0)$ and $(L_1,U_1)$ two slicings of $J$ with $(L_0,U_0)\leq (L_1,U_1)$. Then we have natural isomorphisms
\[
S_{U_0}R_{L_1}\cong R_{L_1}S_{U_0}
\]
Moreover $S_{U_0}R_{L_1}$ is the fiber of the natural transformation $R_{L_1}\to R_{L_0}$ and $R_{L_1}S_{U_0}$ is the cofiber of the natural transformation $S_{U_1}\to S_{U_0}$
\end{lemma}
\begin{proof} 
In the same notation as in the proof of Lemma \refbf{lem.iterated} we have that  $\E_1\subseteq \E_0$ and $\M_0\subseteq \M_1$. Since both $(\E_0,\M_0)$ and $(\E_1,\M_1)$ are normal factorization systems, the isomorphisms of Lemma \refbf{lem.iterated} give the diagram
\[
\begin{kodi}
\obj{
	|(zero)|0 & S_1X &[1cm] S_0X &[1cm] X \\
	&|(zero1)|0& F & R_1X \\
	&& |(zero2)|0 & R_0X & |(zero3)|0\\
};
\mor zero \E_1:-> S_1X {\,\E_0\cap \M_1}:-> S_0X \M_0:-> X \E_1:-> R_1X {\E_0\cap \M_1}:-> R_0X \M_0:-> zero3;
\mor S_1X \E_1:-> zero1 -> F -> zero2 \M_0:-> R_0X;
\mor S_0X -> F -> R_1X;
\end{kodi}
\]
where every square is a pullout. The fact that each class $\E_i$ is closed under pushout and each $\M_i$ is closed under pullback now gives that the arrows $S_0X \to F \to 0$ and $0\to F \to R_1X$ are respectively  the $(\E_1, \M_1)$-factorization of $S_0X\to 0$ and the $(\E_0, \M_0)$-factorization of $0\to R_1X$, so that $R_1 S_0 X \cong F \cong S_0 R_1 X$.

To prove the second part of the statement, notice that by definition of normal factorization system associated to the slicing $(L_0,U_0)$ we have a fiber sequence
\[
\begin{kodi}
\obj{
|(A)| S_0R_1X &[2cm] |(B)| R_1X  \\
|(C)| 0 & |(D)| R_0R_1X  \\
};
\mor A -> B -> D;
\mor * -> C -> *;
\end{kodi}
\]
and the conclusion follows from the natural isomorphism $R_0R_1X\cong R_0X$. Dually one proves the statement on the cofiber of $S_{U_0}\to S_{U_1}$.
\end{proof}

\end{remark}

\subsection{A tale of intervals}
Although a few of the statements we are going to prove hold more generally for arbitrary $\mathbb{Z}$-posets, for the remainder of this section we will restrict our attention to $\mathbb{Z}$-posets which are totally ordered sets.
\begin{definition}
Let $J$ be a poset. An \emph{interval} in $J$ is a subset $I\subseteq J$ such that if $x,y\in I$ and $x\leq z\leq y$ in $J$, then $z\in I$.
\end{definition}

\begin{example}\label{class.is.interval}
Let $J$ be a totally ordered $\mathbb{Z}$-poset, and let $\sim$ be the equivalence relation from Lemma \refbf{equivalence}. For $i\in J$, let $I_i$ be the equivalence class of $i$. Then $I_i$ is an interval. Namely, id $x,y\in I_i$ then there exist integers $a,b$ with $i+a\leq x$ and $y\leq i+b$ so if $x\leq z\leq y$ then $i+a\leq z\leq i_b$ and so $z\sim i$.
\end{example}

Clearly, the intersection of a lower set and an upper set is an interval. Remarkably, in totally ordered sets also the converse is true. Although this is a classical (and easy) result, we recall its proof for completeness.
\begin{lemma}\label{lem.I-is-intersection}
Let $J$ be a totally ordered set. Then a subset $I\subseteq J$ is an interval if and only if $I$ can be written as the intersection of an upper set and a lower set.
\end{lemma}
\begin{proof}
Let 
\[
L_I=\bigcup_{x\in I} (-\infty,x]; \qquad  U_I=\bigcup_{x\in I} [x,+\infty).
\]
Then clearly $L_I$ is a lower set, $U_I$ is an upper set and we have $I\subseteq L_I\cap U_I$. Moreover, if $y\in L_I\cap U_I$ then there exist $x_0,x_1$ in $I$ such that $y\in (-\infty,x_1]\cap [x_0,+\infty)=[x_0,x_1]$. Since $x_0\leq y\leq x_1$ and $I$ is an interval, we have $y\in I$, and so $L_I\cap U_I\subseteq I$.
\end{proof}
\begin{lemma}
In a totally ordered set, the upper set and the lower set intersecting in a nonempty interval $I$ are uniquely determined by $I$.
\end{lemma}
\begin{proof}
Let $I\subseteq J$ be a interval and let 
\[
U_I=\bigcap_{U\supseteq I} U; \qquad L_I=\bigcap_{L\supseteq I} L,
\]
with $U$ and $L$ ranging over the upper sets and the lower sets in $J$ containing $I$, respectively. Then it is clear that $I\subseteq U_I\cap L_I$ and we want to show that actually $I=U_I\cap L_I$ and that if $I=\tilde{U}\cap \tilde{L}$ then $\tilde{U}=U_I$ and $\tilde{L}=L_I$. By Lemma \refbf{lem.I-is-intersection} there exist an upper set $\tilde{U}$ and a lower set $\tilde{L}$ such that $I=\tilde{U}\cap \tilde{L}$. By definition of $U_I$ and $L_I$ we have $U_I\subseteq \tilde{U}$ and $L_I\subseteq \tilde{L}$. Therefore $I\subseteq L_I\cap U_I\subseteq \tilde{L}\cap \tilde{U}=I$ and so $I=U_I\cap L_I$. Now we want to show that $U_I=\tilde{U}$. Since $U_I\subseteq U_0$ we only need to show that $\tilde{U}\subseteq U_I$. Let $x\in \tilde{U}$ and let $y\in I$. Since $J$ is totally ordered, either $x\leq y$ or $x\geq y$. In the first case, since $L_0$ is a lower set, we have $x\in \tilde{L}$ and so $x\in \tilde{L}\cup \tilde{U}=I\subseteq U_I$. In the second case, since $U_I$ is an upper set, we have directly $x\in U_I$.
\end{proof}
By the above lemma, the following definition is well-posed.
\begin{definition}\label{std.endocardium}
Let $J$ be a totally ordered $\mathbb{Z}$-poset and let $\tee\colon \O(J)\to \ts(\C)$ be a $J$-slicing on a stable $\infty$-category $\C$. For every nonempty interval $I=L_I\cap U_I$ in $J$ we set 
\[
\C_I=\C_{L_I}\cap \C_{U_I}.
\]
We also set $\C_\emptyset=\{\mathbf{0}\}$.
\end{definition}
\begin{remark}
The whole of $J$ is an interval, with $L_J=U_J=J$. From \adef\refbf{std.endocardium} we obtain $\C_J=\C$, as expected. Also, every upper set $U$ is an interval, with $U_U=U$ and $L_U=J$. So from  \adef\refbf{std.endocardium} we find that the subcategory of $\C$ associated to $U$ as an interval is precisely the subcategory $\C_U$ associated to $U$ as an upper set. The same happens for lower sets. This shows that the notation introduced in  \adef\refbf{std.endocardium} is consistent with the notation for $J$-slicings.
\end{remark}

\begin{example}\label{Ci}
For every $i,j$ in $J$ with $i\leq j$ one has the four intervals $(i,j)$, $(i,j]$, $[i,j)$, $[i,j]$ and consequently the four subcategories $\C_{(i,j)}$, $\C_{(i,j]}$, $\C_{[i,j)}$ and $\C_{[i,j]}$of $\C$. In particular for every $i\in J$ we have the interval $[i,i]$ consisting of the single element $i$. To avoid cumbersome notation, we will always write $\C_i$ for $\C_{[i,i]}$. The subcategories $\C_i$ with $i$ ranging in $J$ are called the \emph{slices} of the $J$-slicing $\tee$.
\end{example}

\begin{definition}\label{def.bounded}
Let $\tee$ be a $J$-slicing on $\C$. We say that $\C$ is \emph{$J$-bounded} if 
\[
\C=\bigcup_{i,j\in J}\C_{[i,j]}.
\]
Similarly, we say that $\C$ is \emph{$J$-left-bounded} if $\C=\bigcup_{i\in J}\C_{[i,+\infty)}$ and \emph{$J$-right-bounded} if $\C=\bigcup_{i\in J}\C_{(-\infty,i]}$. \end{definition}
\begin{remark}
This notion is well known in the classical as well as in the quasicategorical setting: see \cite{BBDPervers,LurieHA}. In particular, when $\tee$ is a $\Z $-family of $t$-structures on $\C$,  then $\C$ is $\Z $-bounded (resp., $\Z $-left-bounded, $\Z $-right-bounded) if and only if $\C$ is bounded (resp., left-bounded, right-bounded) with respect to the $t$-structure $\tee_0$, agreeing with the classical definition of boundedness as given, \eg, in \cite{BBDPervers}.
\end{remark}
\begin{remark}
Since $\C_{[i,j]}=\C_{[i,+\infty)}\cap \C_{(-\infty,j]}$ one immediately sees that $\C$ is $J$-bounded if and only if $\C$ is both $J$-left- and $J$-right-bounded.
\end{remark}

The following remark is the first step towards the definition of factorization of morphisms associated with interval decompositions of $J$.
\begin{remark}
A nonempty interval in a totally ordered set $J$ is equivalent to the datum of a pair of upper sets $U_0$ and $U_1$ with $U_1\subseteq U_0$, i.e., to the datum of a strictly monotone morphism of posets $\uno \to \O(J)$. Namely, we have seen that $I$ is equivalent to the datum of an upper set $U_I$ and a lower set $L_I$, which are uniquely determined by $I$. Let us set $U_0=U_I$ and $U_1=J\setminus L_I$. Then, since $\O(J)$ is totally ordered by Remark \refbf{rem.totOj}, we have either $U_0\subseteq U_1$ or $U_1\subseteq U_0$. If $U_0\subseteq U_1$ then we have $I=U_0\cap L_I\subseteq U_I\cap L_I=\emptyset$ against the assumption on $I$. So $U_1\subseteq U_0$ and $i\mapsto U_i$ for $i=0,1$ defines a monotone map from $\uno $ to $\O(J)$. Moreover this map is strictly monotone since we have excluded the possibility $U_0\subseteq U_1$ and so we can't have $U_0=U_1$. By removing the assumption that $I$ is nonempty, we can say that an interval in $J$ is given by a (non necessarily strictly monotone) morphism of posets $\uno \to \O(J)$. Actually this is not completely accurate, since all constant maps from $\uno $ to $\O(J)$ will correspond to the empty interval. Yet it will be extremely convenient to always think of intervals as monotone maps to $\O(J)$, so we will systematically adopt this point of view in what follows. In other words we will identify a monotone map $I\colon\uno \to \O(J)$ with the interval $I=U_0\cap L_1$, where $I(0)=(L_0,U_0)$ and $I(1)=(L_1,U_1)$.
\end{remark}

\begin{remark}\label{oi.vs.oj}If $I\colon \uno \to \O(J)$ is an interval in a totally ordered set $J$, then $[U_0,U_1]$ is an interval in the totally ordered set $\mathcal{O}(J)$. It is easy to see that intersecting with $I$ defines a bijection of totally ordered sets
\begin{align*}
[U_0,U_1]&\to \mathcal{O}(I)\\
U&\mapsto U\cap I.
\end{align*}
\end{remark}

\begin{lemma}
Let $I\colon\uno \to \O(J)$ be an interval in $J$, and let $\C_I$ be the corresponding subcategory of $\C$, for a given $J$-slicing. Then the restriction of $S_{U_0}$ to $\C_{L_1}$ and the restriction of $R_{L_1}$ to $\C_{U_0}$ both take values in $\C_I$. 
\end{lemma}
\begin{proof}
We split the proof in two cases. If $I=\emptyset$ then $U_0=U_1$ and so for any $X$ in $\C_{L_1}$ we have $S_0X\cong S_1X=\mathbf{0}$. So $S_{0}\vert_{\C_{L_1}}$ does take its values in $\C_I=\C_\emptyset=\{\mathbf{0}\}$ in this case. If $I\neq \emptyset$, then $U_1\subseteq U_0$. Since $S_{0}$ takes values in $\C_{U_0}$, we only need to show that it maps $\C_{L_1}$ into itself. In other words we want to show that if $X\in \C_{L_1}$ then $S_0X\xrightarrow{\sim} R_1S_0X$. From the fiber sequence 
\[
\begin{kodi}
\obj{
|(A)| S_1S_0X &[1cm] |(B)| S_0X  \\
|(C)| 0 & |(D)| R_1S_0X  \\
};
\mor A -> B -> D;
\mor * -> C -> *;
\end{kodi}
\]
we see we are reduced to showing that $S_1S_0X\cong \mathbf{0}$. Since $U_1\subseteq U_0$, we have $S_1S_0X\cong S_1X$. But, since $X\in \C_{L_1}$ we have $S_1X\cong \mathbf{0}$. This concludes the proof in the case $I\neq \emptyset$. The proof for $R_{L_1}$ is completely analogous.
\end{proof}
By the above lemma and by Lemma \refbf{lem.for-homology} we can give the following
\begin{definition}\label{def.homology}
Let $I\colon\uno \to \O(J)$ be an interval in $J$, and let $\tee\colon \O(J)\to \ts(\C)$ be a $J$-slicing on a stable $\infty$-category $\C$. The functor
\[
\mathcal{H}^I\colon \C\to \C_I
\]
is defined as the composition $\mathcal{H}^I=R_{L_1}S_{U_0}=S_{U_0}R_{L_1}$. 
\end{definition}
As for the functors $R_L$ and $S_U$ we will often implicitly compose $\mathcal{H}^I$ with the inclusion $\C_I\to \C$ and look at it as an endofunctor of $\C$. Notice that if $I$ is the empty interval then $\mathcal{H}^I$ is the zero functor.
\begin{remark}\label{rem.questo}
By looking at a lower set $L$ and to an upper set $U$ as intervals, the above definition gives $\mathcal{H}^L=R_L$ and $\mathcal{H}^U=S_U$. In particular we find
\[
\mathcal{H}^I=\mathcal{H}^{U_0}\mathcal{H}^{L_1}=\mathcal{H}^{L_1}\mathcal{H}^{U_0}
\]
and, by Lemma \refbf{lem.for-homology},  $\mathcal{H}^I$ is the cofiber of the natural transformation $\mathcal{H}^{U_1}\to \mathcal{H}^{U_0}$.
\end{remark}
\begin{remark}
Let $I,\tilde{I}\subseteq J$ two intervals, with $I\subseteq \tilde{I}$. Then
\[
\mathcal{H}^I\mathcal{H}^{\tilde{I}}=\mathcal{H}^{\tilde{I}}\mathcal{H}^I=\mathcal{H}^I.
\] 
Namely, if $I$ is empty, then there is nothing to prove. If $I$ is nonempty, as $I$
is a sub-interval of $\tilde{I}$ we have $U_0\subseteq \tilde{U}_0$ and $L_1\subseteq \tilde{L}_1$. Therefore $(\tilde{L}_0,\tilde{U}_0)\leq (L_0,U_0)\leq (L_1,U_1)\leq (\tilde{L}_1,\tilde{U}_1)$, and so $S_{U_0}R_{\tilde{L}_1}=R_{\tilde{L}_1}S_{U_0}$ by Lemma \refbf{lem.for-homology} as well as $R_{L_1}R_{\tilde{L}_1}=R_{L_1}$ and $S_{U_0}S_{\tilde{U}_0}=S_{U_0}$ by Lemma \refbf{lem.iterated}. Therefore,
\[
\mathcal{H}^I\mathcal{H}^{\tilde{I}}=R_{L_1}S_{U_0}R_{\tilde{L}_1}S_{\tilde{U}_0}=R_{L_1}R_{\tilde{L}_1}S_{U_0}S_{\tilde{U}_0}=R_{L_1}S_{U_0}=\mathcal{H}^I.
\]
The proof that $\mathcal{H}^{\tilde{I}}\mathcal{H}^I=\mathcal{H}^I$ is similar.
\end{remark}

\begin{remark}
If $I$ and $\tilde{I}$ are two disjoint intervals in the totally ordered set $J$ then either every element of $I$ is strictly smaller than every element of $\tilde{I}$ or vice versa. If we are in the first case, then $\C_{I}$ is \emph{right-orthogonal} to $\C_{\tilde{I}}$, \ie, $\C(X,Y)$ is contractible whenever $X\in \C_{\tilde{I}}$ and $Y\in \C_{I}$. Namely, by the assumption on $I$ and $\tilde{I}$ we have $\tilde{U}_0\subseteq U_1$ and so $\C_{\tilde{I}}\subseteq \C_{U_1}$. On the other hand, $\C_I\subseteq \C_{L_1}$ and $\C_{U_1}$ is right-orthogonal to $\C_{L_1}$ by definition of $t$-structure.
\end{remark}
\begin{remark}\label{rem.questaltro}
If $I$ and $\tilde{I}$ are two disjoint intervals in the totally ordered set $J$ , then 
\[
\mathcal{H}^I\mathcal{H}^{\tilde{I}}=\mathcal{H}^{\tilde{I}}\mathcal{H}^I=0.
\] 
Indeed, the statement is trivial if either $I$ or $\tilde{I}$ are empty. When they are nonempty, up to exchanging the role of $I$ and $\tilde{I}$ we may assume that every element of $I$ is strictly smaller than every element of $\tilde{I}$ . Then we have $(L_0,U_0)\leq (L_1,U_1)\leq (\tilde{L}_0,\tilde{U}_0)\leq (\tilde{L}_1,\tilde{U}_1)$ and so
\[
\mathcal{H}^I\mathcal{H}^{\tilde{I}}=R_{L_1}S_{U_0}R_{\tilde{L}_1}S_{\tilde{U}_0}=R_{L_1}R_{\tilde{L}_1}S_{U_0}S_{\tilde{U}_0}=R_{L_1}S_{\tilde{U}_0}=0,
\]
by lemma \refbf{lem.iterated}. Similarly one shows that $\mathcal{H}^{\tilde{I}}\mathcal{H}^I=0$.
\end{remark}
The above Remarks \refbf{rem.questo} and \refbf{rem.questaltro} are actually two particular instances of the following general result. The proof is completely analogous to those in the remarks above, and so it is omitted.
\begin{proposition}\label{H-intersection}
Let $I$ and $\tilde{I}$ be two intervals in a $\mathbb{Z}$-toset, and let $\tee\colon \O(J)\to  \ts(\C)$ be a $J$-slicing on a stable $\infty$-category $\C$. Then
\[
\mathcal{H}^I\mathcal{H}^{\tilde{I}}=\mathcal{H}^{\tilde{I}}\mathcal{H}^I=\mathcal{H}^{I\cap \tilde{I}}.
\]
\end{proposition}
We conclude this Section with a notational convention, which will be useful later. 
\begin{notat}\label{Hi}
Consistently with the notation from Example \ref{Ci},
for every $i$ in $J$ we write $\mathcal{H}^i$ for $\mathcal{H}^{[i,i]}$.
\end{notat}

\section{Interval decompositions and towers.}\label{sec:towers}
\begin{modifyepigraph}{.6}
\epigraph{\cjRL{hAbAh ner:dAh w:nAb:lAh +sAM ,s:pAtAM 'a:+sEr l'o yi+s:m:`U 'iy+s ,s:pat re`ehU;}}{\cjRL{mo+sEh}}
\end{modifyepigraph}
\begin{remark}
For the whole section $(J,\le)$ will be a fixed totally ordered $\Z$-poset and $\tee\colon \O(J)\to  \ts(\C)$ will be a $J$-slicing. 
\end{remark}
\begin{definition}
A \emph{$(k+2)$-fold interval decomposition} of $J$ is a morphism of posets $I_{\ordk}\colon \ordk\to \O(J)$. 
\end{definition}
\begin{notat}
When no ambiguity is possible, the image of $j\in \{0,1,\dots,k\}$ via $I_{\ordk}$ will be denoted simply by $(L_j,U_j)$. For every $j=0,\dots,k+1$ the interval $I_j=U_{j-1}\cap L_j$ is called the \emph{$j$-th interval} in the decomposition, with the convention that $U_{-1}=J=L_{k+1}$. 

The factorization system associated with $(L_j,U_j)$ will be denoted by $(\E_j,\M_j)$. Notice that, since $I_{\ordk}$ is a morphism of posets we have $\E_{j+1}\subseteq \E_j$ and $\M_{j+1}\supseteq \M_j$. 
\end{notat}
This implies that the composition $\ordk \xto{I_{\ordk}} \O(J) \xto{\tee} \ts(\C)$ is a $k$-fold factorization system; in other words
\begin{lemma}\label{k.fold.fact}
Let $(\E_j, \M_j)$ as above. Then every arrow $f\colon X\to Y$ in $\C$ can be uniquely factored into a composition
\[
X \xto{\E_{k}} Z_{k} \xto{\E_{k-1}\cap \M_{k}} Z_{k-1}\to\dots\to Z_{1} \xto{\E_{0}\cap \M_{1}} Z_{0} \xto{\M_{0}} Y.
\]
\end{lemma}
\begin{proof}
Since $(\E_k,\M_k)$ is a factorization system, we have a (unique) factorization $X \xto{\E_{k}} Z_{k} \xto{\M_{k}} Y$.  Since $(\E_{k-1},\M_{k-1})$ is a factorization system, we can (uniquely) factor $Z_k\to Y$ as $Z_k \xto{\E_{k-1}} Z_{k-1} \xto{\M_{k-1}} Y$. Since $\M_{k-1}\subseteq \M_k$, the morphism $Z_{k-1} \to Y$ is also in $\M_k$ and so, by the 3-for-2 property, also $Z_k \to Z_{k-1}$ is in $\M_k$. Therefore $Z_k \xto{\E_{k-1}} Z_{k-1}$ is in $\E_{k-1}\cap \M_k$. then one concludes iterating this argument.
\end{proof}
\begin{definition}
The sequence of morphism in the factorization of $f\colon X\to Y$ in Lemma \refbf{k.fold.fact} is called the \emph{$I_{\ordk}$-tower} of $f$, and it is denoted $\rook(f,I_{\ordk})$, or simply by $\rook(f)$ when the interval decomposition $I_{\ordk}$ is clear from the context.
\end{definition}
\begin{remark}
When $\C$ is the stable $\infty$-category of spectra and $X\to 0$ and $0\to X$ are the terminal and the initial morphism of $X$, respectively, the above notation and construction is in line with the classical Postnikov and Whitehead towers of $X$, i.e., with the sequences
\[
X\to \dots \to R_2X\to R_1X\to R_0X \to 0
\]
\[
0\to \dots \to S_2X\to S_1X\to S_0X \to X
\]
of factorizations obtained from the (stable image of) the $n$-connected factorization system of \cite{Joy}. 
\end{remark}
\begin{remark}\label{rem.trivial-factorizations}
If both $X$ and $Y$ are in $\C_{L_k}$, then the morphism $X\to Z_{k}$ in $\rook(f,I_{\ordk})$ is an isomorphism. Indeed, by construction the morphism $Z_{k}\to Y$ is in $\M_{k}$. Since both $X\to 0$ and $Y\to 0$ are in $\M_{k}$ then also $X\to Y$ is in $\M_{k}$ by 3-for-2, and so also $X\to Z_{k}$ in in $\M_{k}$ again by 3-for-2. But by construction $X\to Z_{k}$ is in $\E_{k}$, so it is an isomorphism. By the same argument one sees that if both $X$ and $Y$ are in $\C_{U_0}$, then the morphism $Z_{0}\to Y$ is an isomorphism.
\end{remark}

\begin{corollary}\label{cor:perPostnikov}
Let $I_{\ordk}$ be a $(k+2)$-fold interval decomposition of $J$. Then for any object $Y$ in $\C$, the tower $\rook\Big(\var{0}{Y},I_{\ordk}\Big)$
is
\[
0 \xto{\E_{k}} \mathcal{H}^{U_k}Y \xto{\E_{k-1}\cap \M_{k}} \mathcal{H}^{U_{k-1}}Y\to\dots\to \mathcal{H}^{U_1}Y \xto{\E_{0}\cap \M_{1}} \mathcal{H}^{U_0}Y \xto{\M_{0}} Y.
\]
Moreover, the arrows $f_j\colon \mathcal{H}^{U_j}Y\to \mathcal{H}^{U_{j-1}}Y$ in $\rook(0\to Y,I_{\ordk})$ are such that $\cofib(f_j)=\mathcal{H}^{I_j}Y\in \C_{I_j}$, for $j=0,\dots,k+1$.
\end{corollary}
\begin{proof}
From the $(k+2)$-fold factorization
\[
0 \xto{\E_{k}} Y_{k} \xto{\E_{k-1}\cap \M_{k}} Y_{k-1}\to\dots\to Y_{1} \xto{\E_{0}\cap \M_{1}} Y_{0} \xto{\M_{0}} Y,
\]
and from the fact that $\E_{0}\supseteq \E_{1} \supseteq\dots\supseteq \E_{k}$ and each class $\E_{j}$ is closed for composition, we see that 
$0\to Y_{j}\to Y$ is the $(\E_j,\M_j)$-factorization of $0\to Y$ and so $Y_{j}=S_{U_j}Y=\mathcal{H}^{U_j}Y$. One concludes by Lemma \refbf{lem.for-homology}. 
\end{proof}
The above corollary motivates the following
\begin{definition}
Let $f\colon X\to Y$ be a morphism in $\C$. A \emph{$I_{\ordk}$-weaved factorization}  for $f$ is a factorization of $f$ of the form
\[
X \xto{f_{k+1}} Z_{k} \xto{f_{k}} Z_{k-1}\to\dots\to Z_{1} \xto{f_{1}} Z_{0} \xto{f_{0}} Y.
\]
with $\cofib(f_j)\in \C_{I_j}$, for $j=0,\dots,k+1$.
\end{definition}
\begin{remark}
If we call $\weaves(f, I_{\ordk})$ the class of $I_{\ordk}$-weaved factorizations for a morphism $f \colon X\to Y$, it is immediate to see that we have a canonical identification
\[
\weaves\Big(\varnobkt{X}{Y},I_{\ordk}\Big) \leftrightarrows \weaves\Big(\varnobkt{0}{\cofib(f)}, I_{\ordk}\Big)
\]
Moreover \acor\refbf{cor:perPostnikov} precisely states that $\rook\Big(\varnobkt{0}{\cofib(f)},I_{\ordk}\Big)$ is an $I_{\ordk}$-weaved factorization, and so we have the existence of $I_{\ordk}$-weaved factorizations for arbitrary morphisms. The remainder of this section is devoted to proving the uniqueness of $I_{\ordk}$-weaved factorizations. This reduces to proving the uniqueness of $I_{\ordk}$-weaved factorizations for initial morphism, i.e., to showing that the only possible $I_{\ordk}$-weaved factorization of $0\to Y$ is its tower.
\end{remark}
\begin{lemma}\label{lemma.vice.versa}
In the above notation, let $f\colon X\to Y$ be a morphism in $\C$. If $X$ is in $\C_{U_j}$ and $\cofib(f)$ is in $\C_{I_j}$ then $0\to X\xto{f}Y$ is the $(\E_j,\M_j)$-factorization of the initial morphism $0\to Y$ and $Y$ is in $\C_{U_{j-1}}$. In particular $f$ is in $\E_{j-1}\cap \M_j$. \end{lemma}
\begin{proof}

Since $X$ is in $\C_{U_j}$, the morphism $0\to X$ is in $\E_j$, and so to show that $0\to X\xto{f} Y$ is the $(\E_j,\M_j)$-factorization of $0\to Y$ we are reduced to showing that $f\colon X\to Y$ is in $\M_j$. Since $\cofib(f)$ is in $\C_{I_j}$, we have in particular that $\cofib(f)\to 0$ is in $\M_j$ and so $0\to \cofib(f)$ is in $\M_j$ by the Sator lemma. Then we have a homotopy pullback diagram
\[
\begin{kodi}
\obj{
	X & 0 \\
	Y &|(cofib)| \cofib(f) \\
};
\mor X :-> 0 {\M_j}:-> cofib <- Y f:<- X;
\end{kodi}
\]
and so $f$ is in $\M_j$ by the fact that $\M_j$ is closed under pullbacks.  To show also that $f\in \E_{j-1}$ let $X\to T\to Y$ be the $(\E_{j-1},\M_{j-1})$-factorization of $f$. Then, since $\E_j\subseteq \E_{j-1}$,  $0\to T\to Y$ is the $(\E_{j-1},\M_{j-1})$-factorization of $0\to Y$. So, by the normality of $(\E_{j-1},\M_{j-1})$ 
we get the diagram
\[
\begin{kodi}
\obj{
0 \\
X &[1cm]|(0')|0 &[1cm]&[1cm]\\
T & U &|(0'')| 0 \\
Y &|(cofib)| \cofib(f) & V &|(0''')| 0 \\
};
\mor[swap] 0 {\E_j}:-> X {\E_{j-1}\cap \M_j}:-> T {\M_{j-1}}:-> Y {\E_j}:-> cofib {\E_{j-1}\cap \M_j}:-> V {\M_{j-1}}:-> 0''';
\mor X {\E_j}:-> 0' {\E_{j-1}\cap \M_j}:-> U {\E_{j-1}\cap \M_j}:-> 0'' {\M_{j-1}}:-> V;
\mor T {\E_j}:-> U {\M_{j-1}}:-> cofib;
\end{kodi}
\]
where all the squares are pullouts, and where we have used the Sator lemma, the fact that $\cofib(f)\to 0$ is in $\M_j$, that the classes $\E$ are closed for pushouts while the classes $\M$ are closed for pullbacks, and the 3-for-2 property for both classes. Since by hypothesis $0\to \cofib(f)$ is in $\E_{j-1}$, we see that $V=0$ and so $T=Y$. Therefore, $Y\in \C_{U_{j-1}}$ and $f\in\E_{j-1}\cap \M_j$.
\end{proof}

\begin{corollary}\label{perPostnikov2} 
Let $Y$ an object in $\C$ and let 
\[
0 \xto{f_{k+1}} Y_{k} \xto{f_k}Y_{{k-1}}\to\dots\to Y_{1} \xto{f_1} Y_{0} \xto{f_0} Y,
\]
be an $I_{\ordk}$-weaved factorization for $0\to Y$. Then $(f_{k+1},\dots,f_1,f_0)=\rook\left(\varnobkt{0}{Y},I_{\ordk}\right)$
\end{corollary}
\begin{proof}
By uniqueness of the $k$-fold factorization we only need to prove that $f_j\in \E_{i_{k-1}}\cap \M_{i_k}$, which is immediate by repeated application of Lemma \refbf{lemma.vice.versa}. 
\end{proof}

\begin{remark}
It's an unavoidable temptation to think of the $I_{\ordk}$-weaved factorization of a morphism $f$ as of its tower $\rook(f,I_{\ordk})$.
As the following counterexample shows, when $f$ is not an initial morphism this is in general not true. Let $J=\Z$, let $k=0$ and let $I_{\ordered{0}}\colon \ordered{0}\to \O(\Z)$ be the slicing $U_0=[0,+\infty)$. Now take a morphism $f\colon X\to Y$ between two elements in $\C_{-1}$. The object $\cofib(f)$ will lie in $\C_{[-1,+\infty)}$, since $\E_0[-1]$ is closed for pushouts, but in general it will not be an element in $\C_{[0,+\infty)}$. In other words, we will have, in general, a nontrivial $(\E_0,\M_0)$-factorization of the initial morphism $0\to \cofib(f)$, i.e., a nontrivial tower $\rook\Big(\varnobkt{0}{\cofib(f)},I_{\ordered{0}}\Big)$. Pulling this back along $Y\to \cofib(f)$ we obtain the $I_{\ordered{0}}$-weaved factorization $X\xto{f_2} Z\xto{f_1} Y$ of $f$, and this factorization will be nontrivial since its pushout is nontrivial. It follows that $(f_2,f_1)$, cannot be $\rook(f,I_{\ordered{0}})$, which is the $(\E_0,\M_0)$-factorization of $f$. Indeed, by the 3-for-2 property of $\M_0$, the morphism $f$ is in $\M_0$, so its $(\E_0,\M_0)$-factorization is trivial. 
\end{remark}

\begin{remark}\label{rem.infittimento1}
Let $I_{\ordered{k'}}\colon \ordered{k'}\to \mathcal{O}(J)$ be a \emph{refinement} of an interval decomposition $I_{\ordk}\colon \ordk\to \mathcal{O}(J)$. This means that for every $i=0,\dots,k+1$ we have an interval decomposition $I_{[\mathbf{k}_i]}\colon [\mathbf{k}_i]\to \mathcal{O}(I_{\ordk;i})$, where $I_{\ordk;i}$ denotes the $i$-th interval in the subdivision $I_{\ordk}$. By the pasting law for pullouts it is immediate to see that for any morphism $f\colon X \to Y$ we have a canonical identification
\[
\weaves\Big(\varnobkt{X}{Y},I_{\ordered{k'}}\Big) \leftrightarrows \weaves\Big(\varnobkt{X}{Y},I_{\ordk}\Big)\times \prod_{i=0}^{k+1}\weaves\Big(\varnobkt{0}{\cofib(f_i)}, I_{[\mathbf{k}_i]}\Big).
\]
\end{remark}

\subsection{Bridgeland slicings} \label{confronto}
\begin{definition}
A $J$-slicing $\tee\colon \O(J)\to  \ts(\C)$ of a stable $\infty$-category $\C$ is called \emph{discrete} if for any object $X$ in $\C$ one has $\mathcal{H}^i(X)=\mathbf{0}$ for every $i$ in $J$ if and only if $X=\mathbf{0}$. A discrete $J$-slicing is said to be of \emph{finite type} if for any object $X$ one has $\mathcal{H}^i(X)\neq\mathbf{0}$ only for finitely many elements $i\in J$.
\end{definition}
\begin{example}\label{example.bounded-t-structure}
A finite type discrete $\mathbb{Z}$-slicing on $\C$ is precisely the datum of a bounded $t$-structure on $\C$.
\end{example}
Suppose now that $\tee$ is of finite type, so that for each $X \in \C$ one has $\mathcal{H}^i(X)=\mathbf{0}$ but for a finite set $\{ i_1^X < \cdots < i_{k_X}^X \} \subseteq J$ of indices $i$, depending on $X$. We can then build up a $(k+2)$-fold interval decomposition $I_{[ \mathbf{k}_X ]}^X$, depending on the object $X$, by setting  $U_j^X=(i_j, +\infty)$. As we are assuming $J$ to be totally ordered, we have $L_j^X=(-\infty, i_j]$. The next proposition shows that the tower of the initial morphism $\mathbf{0} \to X$ associated to this interval decomposition is indeed the ``finest one''. 
\begin{proposition}
Let $\tee\colon \mathcal{O}(J)\to \ts(\C)$ be a $J$-slicing of finite type and let $X$ an object of $\C$. Then for all $j$ we have
\[
\mathcal{H}^{I_j^X}(X)=\mathcal{H}^{i_j^X}(X)
\qquad\text{ and }\qquad  \mathcal{H}^{(i_{j-1}^X,i_{j}^X)}(X)=\mathbf{0}.
\]
\end{proposition}
\begin{proof}
Let us write $i_j$, $I_j$ and  $S_j$ for $i_j^X$, $I_j^X$ and  $S_{U_j}^X$, respectively.
Now, for each $\phi \in J$, using 
By, \aprop\refbf{H-intersection} we have
\[
\mathcal{H}^{i}(\mathcal{H}^{(i_{j-1}, i_{j})}X)=
\begin{cases}
\mathcal{H}^{i}X &\text{if } i\in (i_{j-1},i_{j})\\
\mathbf{0} &\text{otherwise}
\end{cases}.
\]
 But $\mathcal{H}^iX=\mathbf{0}$ for $i \not = i_1 , \cdots , i_{k_X}$, and so $\mathcal{H}^{i}(\mathcal{H}^{(i_{j-1}, i_{j})}X)=\mathbf{0}$ for all $i \in J$. Since $\tee$ is discrete, this gives $\mathcal{H}^{(i_{j-1}, i_{j})}X= \mathbf{0}$, proving the second part of the statement. To prove the first part, recall that
by Lemma \refbf{lem.for-homology}, $\mathcal{H}^{I_j}(X)$ is the cofiber of $S_{j}X\to S_{j-1}X$, while $\mathcal{H}^{i_j}(X)$ is the cofiber of $S_jX\to S_{[i_j,+\infty)}X$. So we are reduced to show that $S_{j-1}X=S_{[i_j,+\infty)}X$. Since $i_{j} > i_{j-1}$ we have $[i_{j},+\infty)\subseteq (i_{j-1},+\infty)=U_{j-1}^X$ and so, by Lemma \refbf{lem.iterated}, $S_{[i_j,+\infty)}X=S_{[i_j,+\infty)}S_{j-1}X$. Thus, we are reduced to show that $S_{j-1}X=S_{[i_j,+\infty)}S_{j-1}X$, i.e., equivalently, that $R_{(-\infty,i_j)}S_{j-1}X=\mathbf{0}$. This is immediate as $
R_{(-\infty,i_j)}S_{j-1}X=\mathcal{H}^{(i_{j-1}, i_j)}X$.
\end{proof}

In particular the above tells us that, writing $\varphi_j$ for $i_{k_X-j}$, the cofiber of the $j$-th morphism of $\rook(\mathbf{0} \to X,I_{[ \mathbf{k}_X ]}^X)$ is $\mathcal{H}^{\varphi_j^X}(X)\in \C_{\varphi_j^X}$. In other words, these towers are weaved factorizations with cofibers in the subcategories $\{ \C_{\varphi} \}_{\varphi \in J}$ and so they correspond to the Harder-Narasimhan filtrations from \cite{Brid}. That is,  Bridgeland's slicings (in their generalized version from \cite{GKR}) are precisely the slicings of finite type in our sense. We show this in detail below.

\begin{definition}\label{def.bridgeland-slicing}
A \textit{Bridgeland $J$-slicing} on $\C$ is a collection $\{\C_{\phi} \}_{\phi \in J}$ of full extension closed sub-$\infty$-subcategories satisfying:  
 \begin{enumerate}[label=$\roman*$)]
\item $\C_{\phi +1}=\C_{\phi}[1]$ for each $\phi \in J$;
\item orthogonality: $\C(X,Y)$ is contractible for each $X \in \C_{\phi}$, $Y \in \C_{\psi}$ for $\phi > \psi$ in $J$;
\item for each object $X \in \C$ there is a finite set $\{ \phi_1 > \cdots > \phi_n \}$ and a factorization of the initial morphism $\mathbf{0} \to X$
$$ \mathbf{0}=X_0 \xrightarrow{\alpha_1} \cdots \xrightarrow{\alpha_n} X_n=X$$
with $\mathbf{0} \neq \cofib(\alpha_i) \in \C_{\phi_i}$ for all $i = 1, \cdots, n$. 
 \end{enumerate}
\end{definition}

\begin{notat}\label{notation.bridg-slicing}
For $\cate{S}$ a subcategory of $\C$, we write $\langle \cate{S}\rangle$ for the smallest extension closed full subcategory of $\C$ containing $\cate{S}$. 
If $M \subseteq J$ is a subset and $\{ \C_{\phi} \}_{\phi \in J}$ is a Bridgeland $J$-slicing on $\C$, we denote $\C_M$ the extension-closed subcategory generated by $\C_{\phi}$ with $\phi \in M$, i.e., we set 
\[
\textstyle \C_M=\left\langle \bigcup_{\phi\in M}\C_\phi \right\rangle.
\]
\end{notat}
\begin{remark}\label{extensions}
Set  $\langle \cate{S}\rangle_0=\mathbf{0}$,  define $\langle \cate{S}\rangle_1$ as the full subcategory of $\C$ generated by $\cate{S}$ and $\mathbf{0}$, and define inductively $\langle \cate{S}\rangle_n$ as the full subcategory of $\C$ on those objects $X$ which fall into a homotopy fiber sequence
\[
\xymatrix{
X_h\ar[r]\ar[d]& X\ar[d]\\
0\ar[r] &X_k
}
\]
with $h,k\geq 1$, $X_h$ in $\langle \cate{S}\rangle_h$, $X_k$ in $\langle \cate{S}\rangle_k$ and $h+k=n$. One clearly has 
\[
\langle \cate{S}\rangle_0\subseteq \langle \cate{S}\rangle_1 \subseteq \langle \cate{S}\rangle_2\subseteq\cdots \subseteq \langle \cate{S}\rangle.
\]
Moreover $\bigcup_n \langle \cate{S}\rangle_n$ is clearly extension closed, so that
\[
\langle \cate{S}\rangle =\bigcup_n \langle \cate{S}\rangle_n.
\] 
\end{remark}
\begin{lemma}\label{closure}
Let $\cate{S}_1,\cate{S}_2$ be two subcategories of $\C$ with $\cate{S}_1\orth \cate{S}_2$, i.e., such that $\C(X,Y)$ is contractible for any $X\in \cate{S}_1$ and any $Y\in\cate{S}_2$. Then $\cate{S}_1\orth \langle\cate{S}_2\rangle$ and $\langle\cate{S}_1\rangle\orth \cate{S}_2$, and so $\langle\cate{S}_1\rangle\orth \langle\cate{S}_2\rangle$
\end{lemma}
\begin{proof}
By Remark \refbf{extensions}, to prove the first statement we are reduced to show that, if $X\in \cate{S}_1$ and $Y\in \langle\cate{S}_2\rangle_n$ then $\C(X,Y)$ is contractible. We prove this by induction on $n$. For $n=0,1$ there is nothing to prove by the assumption $\cate{S}_1\orth \cate{S}_2$. For $n\geq 2$, consider a fiber sequence $Y_h\to Y\to Y_k$ with $1\leq h,k$ and $h+k=n$ as in Remark \refbf{extensions}. Since $\C(X,-)$ preserves homotopy fiber sequences, we get a homotopy fiber sequence of $\infty$-groupoids
\[
\xymatrix{
\C(X,Y_h)\ar[r]\ar[d]& \C(X,Y)\ar[d]\\
{*}\ar[r] &\C(X,Y_k)
}.
\]
By the inductive hypothesis both $\C(X,Y_h)$ and $\C(X,Y_k)$ are contractible, so $\C(X,Y)$ also is. The proof of the second statement is perfectly dual, due to the fact that in $\C$ every fiber sequence is also a cofiber sequence, and $\C(-,Y)$ transforms a cofiber sequence into a fiber sequence.
\end{proof}
\begin{lemma}\label{uno}
Let $(L,U)$ be a slicing of $J$, and let $\C_L$ and $\C_U$ be defined according to Notation \refbf{notation.bridg-slicing}.  Then $\C_U\orth \C_L$.\end{lemma}
\begin{proof}
Since by definition $\C_\phi\orth\C_\psi$ for $\phi>\psi$, the statement immediately follows from Lemma \refbf{closure}.
\end{proof}

\begin{lemma}\label{due}
In the above hypothesis and notation, every object $Y$ of $\C$ sits into a homotopy fiber sequence $Y_{U}\to Y\to Y_{L}$ with $Y_{U}\in \C_{U}$ and $Y_{L}\in \C_{L}$.
\end{lemma}
\begin{proof}
Let
\[
\mathbf{0}=X_0 \xrightarrow{\alpha_1} X_1\cdots \xrightarrow{\alpha_{{\bar\imath}}}X_{{\bar\imath}}\xrightarrow{\alpha_{{\bar\imath}+1}}X_{{\bar\imath}+1}\xrightarrow{}\cdots \xrightarrow{\alpha_n} X_n=X
\]
a factorization of the initial morphism $\mathbf{0} \to X$
with $\mathbf{0} \neq \cofib(\alpha_i) \in \C_{\phi_i}$ for all $i = 1, \cdots, n$, with $\phi_i>\phi_{i+1}$ and with $\phi_{{\bar\imath}}\in U$ and $\phi_{{\bar\imath}+1}\in L$ (with ${\bar\imath}=-1$ or $n$ when all of the $\phi_i$ are in $L$ or in $U$, respectively). 
Consider the pullout diagram
\[
\xymatrix{
X_{{\bar\imath}}\ar[r]\ar[d]_{f_{L}}&0\ar[d]\\
X\ar[r]&\mathrm{cofib}(f_{L})
}\]
together with the 
\[
\mathbf{0}=X_0 \xrightarrow{\alpha_1} X_1\cdots \xrightarrow{\alpha_{{\bar\imath}}}X_{{\bar\imath}} 
\]
and
\[
X_{{\bar\imath}}\xrightarrow{\alpha_{{\bar\imath}+1}}X_{{\bar\imath}+1}\xrightarrow{}\cdots \xrightarrow{\alpha_n} X_n=X.
\]
The first factorization shows that $X_{{\bar\imath}}\in \langle\cup_{i=0}^{{\bar\imath}}\C_{\phi_i}\rangle\subseteq \C_{U}$ while the second factorization shows that $\mathrm{cofib}(f_{L})\in  \langle\cup_{i={\bar\imath}+1}^n\C_{\phi_i}\rangle\subseteq \C_{L}$.
\end{proof}

\begin{lemma}\label{tre}
In the above hypothesis and notation, one has $\C_U[1]\subseteq \C_U$ and $\C_L[-1]\subseteq\C_L$.
\end{lemma}
\begin{proof}
Since the shift functor commutes with pushouts, if an object $X$ is obtained by iterated extensions by objects in $\C_\phi$ with $\phi\in U$, then $X[1]$ is obtained by iterated extensions by objects in $\C_\phi[1]=\C_{\phi+1}$ with $\phi\in U$. In other words, $X[1]$ is an object in $\C_{U+1}$. Since $U$ is an upper set, $U+1\subseteq U$ and so $X[1]\in \C_U$. This proves that $\C_U[1]\subseteq \C_U$. The proof for $\C_L$ is perfectly analogous. 
\end{proof}
The above Lemmas together give the following
\begin{proposition}
In the above hypothesis and notation, the map $\tee\colon (L,U)\mapsto (\C_L,\C_U)$ defines a $J$-slicing of $\C$, i.e., $\tee$ is a $\mathbb{Z}$-equivariant map of posets $\mathcal{O}(J)\to \ts(\C)$. 
\end{proposition}
\begin{proof}
Lemmas \refbf{uno}-\refbf{tre} together precisely say that $(\C_L,\C_U)$ is a $t$-structure on $\C$. Equivariancy of the map is the fact that, as remarked in the proof of Lemma \refbf{tre}, one has $\C_U[1]=\C_{U+1}$. Finally, if $(L_0,U_0)\leq (L_1,U_1)$ then we have $U_1\subseteq U_0$ and so $\C_{U_1}\subseteq \C_{U_0}$, which shows that the map $\tee$ is a morphism of posets.
\end{proof}

\begin{proposition}\label{b-is-discrete}
Let $J$ be a totally ordered $\mathbb{Z}$-poset and let $\C$ be a stable $\infty$-category. Then we have a bijection
\[
\begin{array}{ccc}
 \left\{
\begin{smallmatrix}
\text{finite type discrete}\\
\text{$J$-slicings on $\C$}
\end{smallmatrix}
\right\}
& \longleftrightarrow & 
\left\{
\begin{smallmatrix}
\text{Bridgeland}\\
\text{$J$-slicings on $\C$}
\end{smallmatrix}
\right\}\\[3mm]
(\C_L,\C_U)
& \mapsto & 
\C_{[\phi,+\infty)}\cap \C_{(-\infty,\phi]}\\[3mm]
\Big(\langle \cup_{\phi\in L}\C_\phi \rangle,\langle \cup_{\phi\in U}\C_\phi \rangle\Big)
& \mapsfrom &
\C_\varphi
\end{array}
\]
\end{proposition}
\begin{proof}
The only thing left to be proven is that the above construction actually produces a discrete slicing of finite type. This is actually immediate once one realizes that the factorization  $\mathbf{0}=X_0 \xrightarrow{\alpha_1} \cdots \xrightarrow{\alpha_n} X_n=X$ of the initial morphism $0\to X$ provided by the definition of Bridgeland slicing is actually the weaved factorization corresponding to the interval decomposition of $J$ associated with the decreasing sequence $\phi_1 > \cdots > \phi_n$. One then directly sees that the two constructions indicated in the statement of the proposition are inverse each other.
\end{proof}

\begin{remark}\label{rem:b-is-discrete}
Using the orthogonality condition $\C_\phi\orth\C_\psi$ for $\phi>\psi$,  is not hard to prove by induction on the length of the factorizations that \adef\refbf{def.bridgeland-slicing} actually implies uniqueness of Bridgeland factorizations 
\[\mathbf{0}=X_0 \xrightarrow{\alpha_1} \cdots \xrightarrow{\alpha_n} X_n=X\]
of initial morphisms $0\to X$. As a consequence one has well defined functors
\[
\mathcal{H}^\phi_B\colon \C \to \C_\phi
\]
given by 
\[
\mathcal{H}^\phi_B(X)=\begin{cases}
\cofib(\alpha_i) &\text{ if $\phi=\phi_i$}
\\
\mathbf{0}&\text{ otherwise}
\end{cases}
\]
As it is natural to expect, in the correspondence given by \aprop\refbf{b-is-discrete} the functor $\mathcal{H}^\phi_B$ is identified with the functor $\mathcal{H}^\phi$ associated with the discrete $J$-slicing. Moreover, one easily sees that 
\begin{align*}
\langle \cup_{\phi\in L}\C_\phi \rangle&=\{X\in \C \text{ such that }\mathcal{H}_B^\phi X=\mathbf{0} \text{ for } \phi\in U\}\\
\langle \cup_{\phi\in U}\C_\phi \rangle&=\{X\in \C \text{ such that }\mathcal{H}_B^\phi X=\mathbf{0} \text{ for } \phi\in L\}
\end{align*}
so that the correspondence of \aprop\refbf{b-is-discrete} can be defined sending the pair $(\C_L,\C_U)$ into $\C_{[\phi,+\infty)}\cap \C_{(-\infty,\phi]}$, with inverse
\[
\C_\varphi \mapsto \Big(\{\mathcal{H}_B^\phi X=\mathbf{0} \mid \phi\in U\},\{\mathcal{H}_B^\phi X=\mathbf{0}\mid \phi\in L\}\Big)
\]
\end{remark}

\section{Hearts of $J$-slicings.}\label{hearts}
\begin{modifyepigraph}{.9}
\epigraph{I watched a snail crawl along the edge of a straight razor. That's my dream. That's my nightmare. Crawling, slithering, along the edge of a straight razor\dots\ and surviving.
}{Col. W\@. E\@. Kurtz}
\end{modifyepigraph}

Recall the equivalence relation $x\sim y$ if and only if there are integers $a,b \in \Z$ with $a\leq b$ such that $x + a \le y \le x + b$ on a $\Z$-toset $J$ from lemma \refbf{equivalence}.
\begin{lemma}\label{heartability}
The following are equivalent:
\begin{enumerate}[label=$\roman*$)]
\item the  $\Z$-toset $J$ consists of a single equivalence class with respect to the equivalence relation $\sim$; 
\item there exists an interval $I$ in $J$ such that the map $\varphi\colon (n,x)\mapsto x+n$ is an isomorphism of $\Z$-tosets $\Z\times_{\mathrm{lex}}I\xrightarrow{\sim} J$ (where the $\Z$-action on $I$ is the trivial one and the $\Z$-action on $\Z$ is the translation).
\end{enumerate}
\end{lemma}
\begin{proof}
That $(i)$ implies $(ii)$ is an immediate consequence of \refbf{lemma.representatives}. To prove the converse implication notice that since $\varphi$ is surjective every element in $J$ is equivalent to an element in $I$. So we are reduced to show that all elements in $I$ are equivalent each other. Let $x,y\in I$. We can assume $x\leq y$. Since $J$ is totally ordered we have either $y\leq x+1$ or $x+1\leq y$. In the latter case we have $x\leq x+1\leq y$ and so $x+1\in I$, since $I$ is an interval. But then $\varphi(1,x)=\varphi(0,x+1)$ against the hypothesis on $\varphi$. So we are left with $x\leq y\leq x+1$ which implies $x\sim y$.  
\end{proof}
\begin{definition}
Let $J$ be a $\Z$-toset. A \emph{heart} for $J$ is an interval $J^\heart\subseteq J$ such that $\varphi\colon (n,x)\mapsto x+n$ is an isomorphism of $\Z$-tosets $\Z\times_{\mathrm{lex}}J^\heart \xto{\sim} J$.
\end{definition}
\begin{remark}\label{hasaheart}
Of course, not every $\Z$-toset has a heart. It is easy to see that $J$ has an heart if and only if there is a morphism of $\Z$-tosets $\pi_\heart\colon J \to \Z$. An heart of $J$ is given by $\pi_\heart^{-1}(0)$ in this case.
\end{remark}
\begin{remark}
It is immediate from the definition that $J^\heart$ is a heart of $J$ if and only if $J^\heart+n$ is a heart of $J$, for every $n\in \Z$.
\end{remark}
\begin{example}
If $J=\Z$, with the standard $\Z$-toset structure, the hearts of $J$ are the singletons $\{n\}$ with $n\in \Z$. In particular $\{0\}$ is the standard heart of $\Z$, and all the other hearts are shifts of this. If $J=\mathbb{R}$, with the standard $\Z$-toset structure, then the hearts of $J$ are the intervals of the form $[x,x+1)$ and those of the form $(x,x+1]$, with $x\in \mathbb{R}$.
\end{example}

\begin{example}\label{I.has.a.heart}
Let $(J, \leq)$ be a totally ordered $\mathbb{Z}$-poset, and let $\sim$ be the equivalence relation from Lemma \refbf{equivalence}. For every $i\in J$ let $I_i$ be the equivalence class of $i$. This is an interval in $J$, see Example \refbf{class.is.interval}. Moreover, by Lemma \refbf{lemma.representatives}, $I_i$ has a heart precisely when $i$ is not a fixed point of the $\Z$-action.
\end{example}

\begin{lemma}\label{lemma.plus.one}
Let $I\colon \uno \to \O(J)$ be a heart of $J$. Then $I(1)=I(0)+1$, i.e., $U_1=U_0+1$ (equivalently, $L_1=L_0+1$).
\end{lemma}
\begin{proof}
Assume $U_1\nsubseteq U_0+1$. Then there exists an element $x$ in $U_1\cap(L_0+1)$. Since $I$ is a heart, there exists an element $y$ in $I$ and an integer $n$ such that $x=y+n$. If $n\geq 1$ we have $y+n\in U_0+n\subseteq U_0+1$ and so $x\in (U_0+1)\cap (L_0+1)$, which is impossible. If $n\leq 0$ we have $y+n\in L_1+n\subseteq L_1$ and so $x\in L_1\cap U_1$ which again is impossible. Therefore $U_1\subseteq U_0+1$. Now assume $U_0+1\nsubseteq U_1$. Then there exists an element $x\in (U_0+1)\cap L_1$. Let $y=x-1$. Then $y\in U_0\cap (L_1-1)\subseteq U_0\cap L_1=I$. Since $U_0+1\subseteq U_0$ we also have $x\in I$, and so $\varphi(-1,x)=\varphi(0,y)$, which is impossible. Therefore $U_1= U_0+1$.
\end{proof}
\begin{definition}
Let $J^\heart\subseteq J$ be a heart of $J$ and let $\tee\colon \O(J)\to  \ts(\C)$ be a $J$-slicing on a stable $\infty$-category $\C$. The subcategory $\C_{J^\heart}$ of $\C$ will be called a \emph{heart} of the $J$-slicing $\tee$ and will be denoted $\C^\heart$.
\end{definition}
\begin{notat}
We denote the canonical projection to the heart as $\mathcal{H}^\heart\colon \C\to \C^\heart$; see \adef\refbf{def.homology}.
\end{notat}
\begin{example}
We have seen in Example\refbf{ex.Z-is-t} that a  a $\Z$-slicing on a stable $\infty$-category $\C$ is the same thing as the datum of a $t$-structure $\tee=(\C_{<0},\C_{\geq 0})$ on $\C$. The standard heart $\C^\heart=\C_{\{0\}}$ is called the heart of the $t$-structure $\tee$. The projection to the heart is the functor $\mathcal{H}^0$; see Notation \refbf{Hi}.
\end{example}

From \refbf{cor:perPostnikov} we immediately get the following
\begin{proposition}\label{cor.oss.for.Heart.to.t}
Let $\tee=(\C_{<0},\C_{\geq 0})$ be a bounded $t$-structure on a stable $\infty$-category $\C$, and let $\C^\heart=\C_0$ be its standard heart. Then $\tee$ is completely determined by the functors $\mathcal{H}^j\colon \C \to \C^\heart[j]$ (and so by the functor $\mathcal{H}^0$ alone). More precisely, $\C_{\geq 0}$ is the full subcategory of $\C$ on the objects $Y$ such  that $\mathcal{H}^jY=0$ for any $j< 0$, while $\C_{< 0}$ is the full subcategory of $\C$ on those objects $Y$ for which $\mathcal{H}^jY=0$ for any $j\geq 0$.
\end{proposition}

\begin{remark}\label{evocative}
There is a rather evocative pictorial representation of the heart $\C_{[0,1)}$ of an $\mathbb{R}$-slicing, manifestly inspired by \cite{Brid}: if we depict $\C_{<0}$ and $\C_{\geq 0}$ as contiguous half-planes (refer to Figure \refbf{figure:slices})
then the action of the shift functor can be represented as an horizontal shift, and the closure properties of the two classes $\C_{<0}, \C_{\geq 0}$ under positive and negative shifts are a direct consequence of the shape of these areas. With these notations, an object $Z$ is in the heart $\C_{[0,1)}$ if it lies in a ``boundary region'', \ie if it lies in $\C_{\geq 0}$, but $Z[-1]$ lies in $\C_{<0}$.
\begin{center}
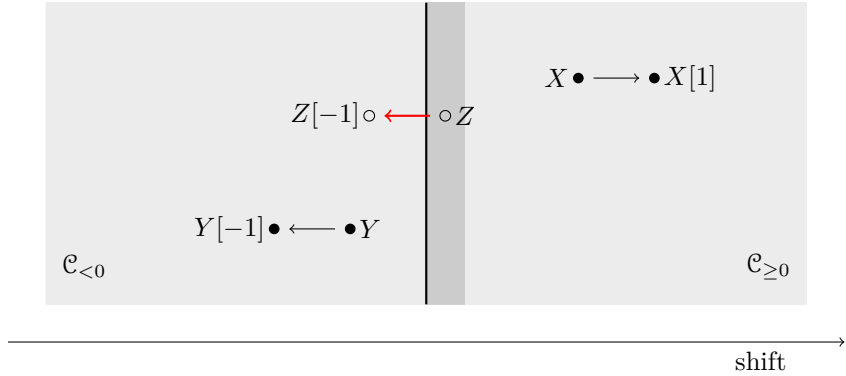
\begin{figure}[h]
\begin{tikzpicture}
\filldraw[gray!15] (5,-2) -- (-5,-2) -- (-5,2) -- (5, 2) -- cycle;
\filldraw[gray!40] (.5, -2) -- (0,-2) -- (0,2) -- (.5,2) -- cycle;
\draw[thick] (0,-2) -- (0,2);
\fill (2,1) circle (2pt) node[left] (X) {$X$};
\fill (-1,-1) circle (2pt) node[right] (Y) {$Y$};
\node at (4.5,-1.5) {$\C_{\geq 0}$};
\node at (-4.5,-1.5) {$\C_{<0}$};
\draw[->] (2.2,1) -- (2.8,1);
\draw[->] (-1.2,-1) -- (-1.8,-1);
\fill[xshift=1cm] (2,1) circle (2pt) node[right] (X) {$X[1]$};
\fill[xshift=-1cm] (-1,-1) circle (2pt) node[left] (Y) {$Y[-1]$};
\draw (.25,0.5) circle (2pt) node[right] {$Z$};
\draw (-.75,0.5) circle (2pt) node[left] {$Z[-1]$};
\draw[thick, ->, red] (.05,0.5) to (-.55,0.5);
\draw[->, xshift=-.5cm, yshift=-.5cm] (-5,-2) -- (6,-2) node[below, pos=.9] {\text{shift}};
\end{tikzpicture}
\caption{Heart of an $\mathbb{R}$-slicing}
\label{figure:slices}
\end{figure}
\end{center}
\end{remark}

Let now $J^\heart\subseteq J$ be a heart, and let $\C^\heart=\C_I$ be the corresponding subcategory of $C$, relative to a given $J$-slicing $\tee$. Writing $I$ as $I\colon \uno \to \O(J)$, for any $n\in \Z$ and any $k\geq 0$ we can consider the interval decomposition $I_{\ordered{k}}\colon \ordered{k}\to \O(J)$ defined by setting $I_{\ordered{k}}(j)=I(0)+j+n$ for $0\leq j\leq k$. By Lemma \refbf{lemma.plus.one}, this corresponds to the collection of $n$ contiguous intervals $I+n,I+n+1,\dots,I+n+k-1\subseteq J$. The corresponding subcategories of $\C$ will be $\C^\heart[n]$, $\C^\heart[n+1]$,\dots,$\C^\heart[n+k-1]$.

 The existence and uniqueness of $I_{\ordered{k}}$-weaved factorizations then specializes to the following statement
\begin{proposition}\label{prop.first-tower}
Let $f\colon X\to Y$ be a morphism in $\C$. Then for any integer $n$ and any positive integer $k$ there exists a unique factorization
\[
X \xto{f_{n+k}} Z_{n+k-1} \xto{f_{n+k-1}}Z_{n+k-2}\to\dots\to Z_{n+1} \xto{f_{n}} Z_{n} \xto{f_{n-1}} Y,
\]
of $f$ such that
$\cofib(f_j)\in \C^\heart[j]$ for any $j=n,\dots,n+k-1$,  $\cofib(f_{n-1})\in \C_{L_0}[n]$  and $\cofib(f_{n+k})\in \C_{U_k}[n]=\C_{U_0}[k+n]$.
\end{proposition}
The content of \aprop\refbf{prop.first-tower} becomes more interesting when $\C$ is \emph{bounded} with respect to the $J$-slicing $\tee$ (see \adef\refbf{def.bounded}), as in this case $\cofib(f)$ lies in $\C_{U_0}[n]$ for $n\ll 0$ and in $\C_{L_0}[n]$ for $n\gg0$ . Namely, if $\C$ is $J$-bounded, then $\cofib(f)$ lies in $\C_{\geq i}$ for some $i\in J$. Since $J^\heart$ is a heart, there exists an element $x\in I$ and an integer $n_0$ such that $i=x+n_0$, so that $\cofib(f)\in \C_{\geq x}[n_0]$. As $x\in I$ we have $x\in U_0$ and so $[x,+\infty)\subseteq U_0$ therefore $\cofib(f)\in \C_{U_0}[n_0]$ and so $\cofib(f)\in \C_{U_0}[n]$ for any $n\leq n_0$. Dually one proves the statement for $n\gg0$. 
 As an immediate consequence, by Remark \refbf{rem.trivial-factorizations} we see that in the $J$-bounded case
if $f\colon X\to Y$ is any morphism in $\C$, then  
 the morphisms $0 \to \cofib(f)_{n+k-1}$ and $\cofib(f)_{n} \to \cofib(f)$ in $\rook\Big(\varnobkt{0}{\cofib(f)},I_{\ordered{k}}\Big)$ are isomorphisms for $n\ll 0$ and $k \gg 0$. Since $\rook\Big(\varnobkt{0}{\cofib(f)},I_{\ordered{k}}\Big)$ is the $I_{\ordered{k}}$-weaved factorization of $0\to \cofib(f)$ by Corollary \refbf{cor:perPostnikov}, and since isomorphisms are preserved by pullouts we see that both $X \xto{f_{n+k}} Z_{n+k-1} $ and $Z_{n} \xto{f_{n-1}} Y$ are isomorphisms.
 
 \par
 This leads to the following
\begin{proposition}\label{prop.Z.Postnikov}
Let $\C$ be a stable $\infty$-category which is bounded with respect to a given $J$-slicing $\tee$. Let $J^\heart\subseteq J$ be a heart for $J$ and let $\C^\heart$ be the corresponding heart in $\C$.
Then for any morphism $f\colon X\to Y$  in $\C$ there exists an integer $n_0$ and a positive integer $k_0$ such that for any integer $n\leq n_0$ and any positive integer $k$ with $k\geq n_0-n+k_0$ there exists a unique factorization of $f$ 
\[
X \xto{\sim} Z_{n+k-1} \xto{f_{n+k-1}}Z_{n+k-2}\to\dots\to Z_{n+1} \xto{f_{n}} Z_{n} \xto{\sim} Y
\]
such that
$\cofib(f_j)\in \C^\heart[j]$ for any $j=n,\dots,n+k-1$.
\end{proposition}
\begin{remark}\label{oss.Z.Postnikov}
By uniqueness in \aprop\refbf{prop.Z.Postnikov}, one has a well defined $\Z $-factorization
\[
X=\lim(Z_j) \to\cdots \to Z_{j+1} \xto{f_{j}} Z_{j} \xto{f_{j-1}}Z_{j-1}\to \cdots\to \mathrm{colim}(Z_j)=Y
\]
with $j$ ranging over the integers, $\cofib(f_j)\in \C^\heart[j]$ for any $j\in \Z $ and with $f_m$ being an isomorphism for $|j|\gg 0$. We will refer to this factorization as the $\C^\heart$-weaved factorization of $f$. Notice how the boundedness of $\C$ has played an essential role: when $\C$ is not bounded, one still has towers for interval decomposition $I_{\ordered{k}}\colon j\mapsto I(0)+j+n$ for arbitrary $k$ and $n$, but in general they do not stabilize.
\end{remark}
\begin{remark}\label{oss.for.Heart.to.t}
Let $J^\heart\subseteq J$ be a heart, and let $(L,U)$ be the slicing $J^\heart(0)$ of $J$ and $\tee=(\C_U,\C_L)$ be the corresponding $t$-structure on $\C$. By Lemma \refbf{lemma.plus.one} the standard heart of $\tee$ is precisely $\C^\heart$. Moreover,
by Corollaries \refbf{cor:perPostnikov} and \refbf{perPostnikov2}, the $\C^\heart$-weaved factorization 
\[
0=\lim(Y_j) \to\cdots \to Y_{j+1} \xto{f_{j}} Y_{j} \xto{f_{j-1}}Y_{j-1}\to \cdots\to \mathrm{colim}(Y_j)=Y
\]
of an initial morphism $0\to Y$ is such that $\cofib(f_j)=\mathcal{H}^jY$ for any $j\in \Z$. Therefore we see that the $t$-structure $(\C_L,\C_U)$ can be completely read from $\C^\heart$-weaved factorizations: an object $Y$ is in $\C_{U}$ if and only if the $\C^\heart$-weaved factorization of $0\to Y$
satisfies $\cofib(f_j)=0$ for any $j< 0$, while $Y$ is in $\C_{L}$ if and only if  $\cofib(f_j)=0$ for any $j\geq 0$. We are going to use this fact later to characterize hearts of $t$-structures among full subcategories of $\C$.
\end{remark}

\begin{remark}\label{bounded.is.bounded}
The content of Remark \refbf{oss.for.Heart.to.t} can be elegantly expressed in terms of the functoriality of slicings (Remark \refbf{rem.slicing-functor}). Namely, by Remark \refbf{hasaheart}, to give a heart $J^\heart$ of $J$ is equivalent to giving  a $\Z$-equivariant morphism of $\pi_\heart\colon \Z$-tosets $J\to\Z$, and by functoriality this induces a map
\[
\pi_\heart\colon \slicings(J,\C)\to \ts(\C).
\]   
The heart $\C^\heart$ of a $J$-slicing $\tee$ on $\C$ is then the standard heart of the corresponding $t$-structure $\pi_\heart(\tee)$. Moreover, one easily sees that $\tee$ is a bounded $J$-slicing if and only if $\pi_\heart(\tee)$ is a bounded $t$-structure. Indeed, since $\pi_\heart$ is a morphism of tosets, for any $i_1\leq i_2$ in $J$ we have $[i_1,i_2]\subseteq \pi_\heart^{-1}[\pi_\heart(i_1),\pi_\heart(i_2)]$. Vice versa, since  $\pi_\heart$ is $\Z$-equivariant, for every $j\in J$ we have $\pi_\heart(j+n)=\pi_\heart(j)+n$ and so for every $n_1\leq n_2$ in $\Z$ there exist $j_1\leq j_2$ in $J$ such that $\pi_\heart(j_1)< n_1$ and $\pi_\heart(j_2)> n_2$. Let $j\in J$ be such $n_1\leq \pi_\heart(j)\leq n_2$. If $j>j_2$ then we would have $\pi_\heart(j)\geq \pi_\heart(j_2)>n_2$, and so $j\leq j_2$. Similarly $j\geq j_1$ and so $\pi_{\heart}^{-1} [n_1,n_2]\subseteq [j_1,j_2]$.
\end{remark}

\subsection{Abelianity of the heart.}
In the following section we present a proof of the fact that a heart $\C^\heart$ of a $J$-slicing on a stable $\infty$-category $\C$
is an abelian $\infty$-category. 

In other words, $\C^\heart$ is homotopy equivalent to its homotopy category $h\C^\heart$, which is an abelian category; this is the higher-categorical counterpart of a classical result, first proved in \cite[\athm\textbf{1.3.6}]{BBDPervers}, 
which 
only relies on properties stated in terms of normal torsion theories in a stable $\infty$-category. 

We begin with the following
\begin{definition}[Additive $\infty$-category]\label{df:additive}
An \emph{additive $\infty$-category}  is an additive
category regarded as an $\infty$-category.
\end{definition}
\begin{remark}
Equivalently, an additive $\infty$-category is a quasicategory $\cate{A}$
such that
\begin{enumerate}[label=$\roman*$)]
\item the hom space $\cate{A}(X,Y)$ is a homotopically discrete infinite loop space for any $X, Y$, \ie, there exists an infinite sequence of $\infty$-groupoids $\Xi_0, \Xi_1,\Xi_2,\dots$, with $\Xi_0\cong \C(X,Y)$ and homotopy equivalences $\Xi_i\cong \Omega \Xi_{i+1}$ for any $i\geq 0$, such that $\pi_n \Xi_0=0$ for any $n\geq 1$;
\item $\cate{A}$ has a zero object, and (homotopy) biproducts.
\end{enumerate}
\end{remark}
\begin{remark}\label{orthogonality-of-a}
If $\cA$ is a full additive $\infty$-subcategory of a stable $\infty$-category $\C$, then $\cA[n_1]\orth\cA[n_2]$ for any $n_1>n_2$. Namely, for an additive $\cA$, let $X\in \cA[n_1]$ and $Y\in\cA[n_2]$. Then $X=Z_1[n_1]$ and $Y=Z_2[n_2]$ for suitable $Z_1,Z_2\in \cA$ and so 
\begin{align*}
\C(X,Y)&=\C(Z_1[n_1],Z_2[n_2])\cong  \C(Z_1,Z_2[n_2-n_1])\\
&\cong \Omega^{n_1-n_2}\C(Z_1,Z_2)=\Omega^{n_1-n_2}\cA(Z_1,Z_2),
\end{align*}
where in the last equality we used the fact that $\cA$ is full. Since $n_1-n_2>0$, the space $\Omega^{n_1-n_2}\cA(Z_1,Z_2)$ is contractible by definition of additive $\infty$-category.
\end{remark}
\begin{remark}
In the triangulated setting one does not have the compatibility between the looping operation on objects and on spaces of morphisms $\Omega\C(X,Y)=\C(X,\Omega Y)$ and so the othogonality condition $\cA[n_1]\orth\cA[n_2]$ for $n_1>n_2$, when needed, has to be imposed by hands. This is nothing but the usual vanishing of negative Ext's one frequently meets as a property of `good' additive subcategories of a triangulated category.  
\end{remark}

\begin{definition}[Abelian $\infty$-category]\label{df:abelinfty}
An \emph{abelian $\infty$-category}  is an abelian
category regarded as an $\infty$-category.
\end{definition}
\begin{remark}
Equivalently, an abelian $\infty$-category is an additive $\infty$-category $\cate{A}$
such that
\begin{enumerate}[label=$\roman*$)]
\item $\cate{A}$ has (homotopy) kernels and cokernels;
\item for any morphism $f$ in $\cate{A}$, the natural morphism from the \emph{coimage} of $f$ to the \emph{image} (see \adef\refbf{imcoim}) of $f$ is an equivalence.
\end{enumerate}
\end{remark}

The rest of the section is devoted to the proof of the following result:
\begin{theorem}\label{heart.is.abelian}
Let $\C^\heart$ be a heart of a $J$-slicing $\tee$ on a stable $\infty$-category $\C$. Then $\C^\heart$ is an abelian $\infty$-category.\end{theorem}
\begin{remark}
If $J=\Z$, so that $\tee$ is the datum of a $t$-structure on $\C$, and $\C^\heart=\C_{\{0\}}$ is the standard heart of $\C$, the the  homotopy category $h\C^\heart$ is the abelian category arising as the standard heart of the $t$-structure $h(\tee)$ on the triangulated category $h\C$.
\end{remark}
In what follows, let $I\colon \uno \to \O(J)$ be an interval such that $\C^\heart=\C_I$. The two factorization systems associated with $I$ will be denoted by $(\E_0,\M_0)$ and $(\E_1,\M_1)$, respectively. By Lemma \refbf{lemma.plus.one} we have $\M_1=\M_0[1]$ and $\E_1=\E_0[1]$.
\begin{lemma}
For any $X$ and $Y$ in $\C^\heart$, the hom space $\C^\heart(X,Y)$ is a homotopically discrete infinite loop space.
\end{lemma}
\begin{proof}
Since $\C^\heart$ is a full subcategory of $\C$, we have $\C^\heart(X,Y)=\C(X,Y)$, which is an infinite loop space since $\C$ is a stable $\infty$-category. 

So we are left to prove that $\pi_n\C(X,Y)=0$ for $n\geq 1$. Since $\pi_n\C(X,Y)=\pi_{n-1}\Omega\C(X,Y)=\pi_{n-1}\C(X,Y[-1])$, this is equivalent to showing that 
$\C(X,Y[-1])$ is contractible. Since $X$ and $Y$ are objects in $\C^\heart$, we have $X\in \C_{U_0}$ and $Y[-1]\in \C_{L_1[-1]}=\C_{L_1-1}=\C_{L_0}$. But $\C_{L_0}$ is right object-orthogonal to $\C_{U_0}$, therefore $\C(X,Y[-1])$ is contractible.
\end{proof}

The subcategory $\C^\heart$ inherits the $0$ object and biproducts (in fact, all finite limits) from $\C$, so in order to prove it is is abelian we are left to prove that it has kernels and cokernels, and that the canonical morphism from the coimage to the image is an equivalence.
\begin{lemma}\label{lemma.qua.e.la}
Let $f\colon X\to Y$ be a morphism in $\C^\heart$. Then $\mathrm{fib}(f)$ is in $\C_{L_1}$ and $\cofib(f)$ is in $\C_{U_0}$.
\end{lemma}
\begin{proof}
Since both $X\to 0$ and $Y\to 0$ are in $\M_1$, by the 3-for-2 property also $f$ is in $\M_1$. Since $\M_1$ is closed for pullbacks, $\mathrm{fib}(f)\to 0$ is in $\M_1$ and so $\mathrm{fib}(f)$ is in $\C_{L_1}$. The proof for $\cofib(f)$ is completely dual.
\end{proof}
\begin{definition}
Denote by
\[
\xymatrix{
0\ar[r]^-{\E_0}&\ker(f)\ar[r]^{\M_0}&\fib(f)
}
\]
the $(\E_0,\M_0)$-factorization of the morphism $0\to \fib(f)$
and by
\[
\xymatrix{
\cofib(f)\ar[r]^{\E_1}&\coker(f)\ar[r]^-{\M_1} &0
}
\]
the $(\E_1,\M_1)$-factorization of the morphism $\cofib(f)\to 0$. We call $\ker(f)$ and $\coker(f)$ respectively the \emph{kernel} and the \emph{cokernel} of $f$ in $\C^\heart$.
\end{definition}
\begin{remark}\label{oss.miracle}
Since $\cofib(f)[-1]\cong \fib{f}$ and $(\M_1,\E_1)[-1]=(\M_0,\E_0)$, one can equivalently define $\coker(f)$ by declaring the $(\E_0,\M_0)$-factorization of $\fib(f)\to 0$ to be $\fib(f)\xto{\E_0}\coker(f)[-1]\xto{\M_0} 0$. Similarly, one can define $\ker(f)$ by declaring the $(\E_1,\M_1)$-factorization of $0\to \cofib(f)$ to be $0\xto{\E_1}\ker(f)[1]\xto{\M_1} \cofib(f)$.
By normality of the factorization system we therefore have the homotopy commutative diagram 
\[
\begin{kodi}
\obj{
	0  &|(ker)| \ker(f) &[1cm] |(fib)| \fib(f)\\
	&|(0')| 0 &|(coker)| \coker(f)[-1] & |(0'')| 0 \\
};
\mor 0 {\E_0}:-> ker \M_0:-> fib \E_0:-> coker \M_0:-> 0'';
\mor ker \E_0:-> 0' \M_0:-> coker \M_0:-> 0'';
\end{kodi}
\]
whose square sub-diagram is a homotopy pullout.
\end{remark}
\begin{lemma}
Both $\ker(f)$ and $\coker(f)$ are in $\C^\heart$.
\end{lemma}
\begin{proof}
By construction $\ker(f)$ is in $\C_{U_0}$, so we only need to show that $\ker(f)$ is in $\C_{L_1}$. By definition of $\ker(f)$, we have that $\ker(f)\to \fib(f)$ is in $\M_0$. Since $\M_0[-1]\subseteq \M_0$, we have that also $\ker(f)[-1]\to \fib(f)[-1]$ is in $\M_0$.
By Lemma \refbf{lemma.qua.e.la}, $\fib(f)[-1]\to 0$ is in $\M_1[-1]=\M_0$ and so we find that also $\ker(f)[-1]\to 0$ is in $\M_0$. 
The proof for $\coker(f)$ is perfectly dual.
\end{proof}
By definition of $\ker(f)$ and $\coker(f)$, the defining diagram of $\fib(f)$ and $\cofib(f)$ can be enlarged as
\[
\begin{kodi}
\obj{
0 & |(ker)| \ker(f) & |(fib)| \fib(f) & X & |(01)| 0 \\
& &|(02)| 0 & Y & |(cofib)| \cofib(f) & |(coker)| \coker(f) & |(03)| 0\\	
};
\mor 0 -> ker -> fib -> X -> 01 -> cofib -> coker -> 03;
\mor fib -> 02 -> Y -> cofib;
\mor X f:-> Y;
\mor ker k_f:{bend left=40},-> X; \mor Y c_f:swap,{bend right=40},-> coker;
\end{kodi}
\]
where $k_f$ and $c_f$ are morphisms in $\C^\heart$.
\begin{definition}\label{imcoim}
Let $f\colon X\to Y$ be a morphism in $\C^\heart$. The \emph{image} $\im(f)$ and the \emph{coimage} $\coim(f)$ of $f$ are defined as $\im(f)=\ker(c_f)$ and  $\coim(f)=\coker(k_f)$. 
\end{definition}
The following lemma shows that $\ker(f)$ does indeed have the defining property of a kernel:
\begin{lemma}\label{is.a.kernel}
The homotopy commutative diagram
\[
\xymatrix{
\ker(f)\ar[r]^{k_f}\ar[d]&X\ar[d]^{f}\\
0\ar[r]&Y}
\]
is a pullback diagram in $\C^\heart$.
\end{lemma}
\begin{proof}
A homotopy commutative diagram 
\[
\xymatrix{
K\ar[r]\ar[d]&X\ar[d]^{f}\\
0\ar[r]&Y}
\]
between objects in the heart is in particular a homotopy commutative diagram in $\C$ so it is equivalent to the datum of a morphism $k'\colon K\to \fib(f)$ in $\C$, with $K$ an object in $\C^\heart$. By the orthogonality of $(\E_0,\M_0)$, this is equivalent to a morphism $\tilde{k}\colon K\to\ker(f)$:
\[\xymatrix{
0\ar[r]\ar[d]_{\E_0}&\ker(f)\ar[d]^{\M_0}\\
K\ar[ru]^{\tilde{k}}\ar[r]_{k'}&\fib(f)
}. \qedhere 
\]
\end{proof}
There is, obviously, a dual result showing that $\coker(f)$ is indeed a cokernel.
\begin{lemma}
The homotopy commutative diagram
\[\xymatrix{
X\ar[r]\ar[d]_{f}&0\ar[d]\\
Y\ar[r]^-{c_f}&\coker(f)
}\]
is a pushout diagram in $\C^\heart$.
\end{lemma}

\begin{proposition}\label{pullout.is.pullout}
A homotopy commutative diagram
\[
\xymatrix{
X\ar[r]^{f}\ar[d]&Y\ar[d]_{g}\\
\zero\ar[r]&Z
}
\]
with vertices in $\C^\heart$ is a pullout diagram in $\C^\heart$ if and only if it is a pullout diagram in $\C$. 
\end{proposition}
\begin{proof}
Since $\C^\heart$ is a full subcategory of $\C$ it is clear that if the given diagram is a pullout in $\C$ then it is automatically also a pullout in $\C^\heart$. Conversely, assume the given diagram is a pullout in $\C^\heart$. This means that $X=\ker(g)$ and $Z=\coker(f)$. By definition of $\ker$ and $\coker$ this implies that we have the following homotopy commutative diagram in $\C$ where each square is a pullout (in $\C$):
\[
\begin{kodi}
\obj{
   |(01)|\zero &[1cm]  X &[1cm]  \fib(g)  &[1cm] Y \\
&  |(02)|\zero &  W & \cofib(f)     \\
&& |(03)|\zero &  Z                 \\
};
\mor 01 \E_0:-> X \E_0:-> 02 -> W \E_0:-> 03 -> Z;
\mor X \M_0:-> {fib g} \E_0:-> W -> {cofib f} \E_1:-> Z;
\mor {fib g} -> Y -> {cofib f};
\mor X f:dashed,{bend left},-> Y g:dashed,{bend left=50},-> Z;
\end{kodi}
\]
Here $X\to \zero$ is in $\E_0$ by definition of $\ker(g)$ and by the Sator lemma, and so also $\fib(g)\to W$ is in $\E_0$ since $\E_0$ is closed under pushouts.
The morphism $\cofib(f)\to Z$ is in $\E_1$ by definition of $\coker(f)$. From the pullout diagram in $\C$
\[
\xymatrix{
\cofib(f)\ar[r]\ar[d] & \zero\ar[d]\\
Z\ar[r] & W[1]
}
\]
 and the fact $\E_1$ is closed under pushouts that follows that $\zero\to W[1]$ is in $\E_1$. Therefore $\zero\to W$ is in $\E_1[-1]=\E_0$ and so $W\to \zero$ is in $\E_0$. This shows that $\fib(g)\to \zero$ is in $\E_0$ and so $X\to \fib(g)$ is an isomorphism. Therefore the given diagram is a pullback in $\C$ and so a pullout in $\C$. 
\end{proof}

\begin{lemma}\label{lemma.titanic}
For $f\colon X\to Y$ a morphism in $\C$, there is a homotopy commutative diagram where all squares are homotopy pullouts:
\[
\begin{kodi}
\obj{
|(ker)| \ker(f) &[1cm] |(fib)| \fib(f) &[1cm] X &[1cm] 0 &[1cm]\\
|(01)| 0 &|(coker-1)| \coker(f)[-1] & |(Z)| Z_f &|(ker1)| \ker(f)[1] &|(02)| 0 \\
& |(03)| 0 & Y & |(cofib)| \cofib(f) & |(coker)| \coker(f) \\
};
\mor ker -> fib -> X -> 0 {{\E_1}}:-> ker1 -> 02 {{\M_1}}:-> coker;
\mor ker {\E_0}:-> 01 -> coker-1 -> Z -> ker1 {{\M_1}}:-> cofib -> coker;
\mor fib {\E_0}:-> coker-1 {\M_0}:-> 03 -> Y -> cofib;
\mor[swap] X {\E_0}:-> Z {{\M_1}}:-> Y;
\mor X f:{near start},dashed,{bend left},-> Y; 
\mor ker k_f:{bend left},-> X; 
\mor[swap] Y c_F:{bend right},-> coker;
\end{kodi}\]
uniquely determining an object $Z_f\in \C^\heart$.
\end{lemma}
\begin{proof}
Define $Z_f$ as the homotopy pullout
\[
\xymatrix{
\fib(f)\ar[r]\ar[d]_{\E_0} \pp & X\ar[d]^{\E_0}\\
\coker(f)[-1]\ar[r]& Z_f}
\]
Here the vertical arrow on the right is in $\E_0$ since the vertical arrow on the left is in $\E_0$ by definition of $\coker(f)$ (see Remark \refbf{oss.miracle}) and $\E_0$ is preserved by pushouts. Next, paste on the left of this diagram the pullout given by Remark \refbf{oss.miracle} and build the rest of the 
diagram by taking pullbacks or pushouts, recalling that $(\E_0,\M_0)[1]=(\E_1,\M_1)$. Use again Remark \refbf{oss.miracle} and the fact that $\M_1$ is closed under pullbacks to see that $Z_f\to Y$ is in $\M_1$. 
Finally, we have 
\[
0\xto{\E_0} X\xto{\E_0} Z_f\xto{\M_1} Y\xto{\M_1} 0,
\]
and so $Z_f$ is in $\C^\heart$.
\end{proof}
\begin{proposition}\label{im.iso.coim}
There is an isomorphism $\im (f)\cong\coim (f)$.\end{proposition}
\begin{proof}
By definition, $\im(f)$ and $\coim(f)$ are defined by the factorizations
\[
\xymatrix{
0\ar[r]^-{\E_0}&\im(f)\ar[r]^{\M_0}&\fib(c_f)
}
\]
and
\[
\xymatrix{
\cofib(k_f)\ar[r]^{\E_1}&\coim(f)\ar[r]^-{\M_1} &0
}
\]
The diagram in Lemma \refbf{lemma.titanic} shows that we have $\fib(c_f)=Z_f=\cofib(k_f)$. Therefore, 
what we need to exhibit are the $(\E_0,\M_0)$ factorizations of $0\to Z_f$ and the $(\E_1,\M_1)$ factorization of $Z_f\to 0$. Since $Z_f$ is an object in $\C^\heart$, these are
\[
0\xrightarrow{\E_0}Z_f\xrightarrow{\mathrm{id}_{Z_f}}Z_f
\]
and 
\[
Z_f\xrightarrow{\mathrm{id}_{Z_f}}Z_f\xrightarrow{\M_1}0,
\]
respectively, thus giving $\im(f)\cong Z_f\cong \coim(f)$.
\end{proof}

\subsection{Abelian subcategories as hearts.}
By the results in the previous section and by Corollary \refbf{cor.oss.for.Heart.to.t} we see that hearts of bounded $t$-structures are very peculiar subcategories of a stable $\infty$-category $\C$: they are abelian and every morphism in $\C$ admits a unique $\C^\heart$-weaved factorization. As we are going to show, these two properties precisely characterize hearts among full subcategories of $\C$.
\begin{definition}
Let $f\colon X\to Y$ be a morphism in $\C$, and let $\cA$ be an additive $\infty$-subcategory of $\C$. An $\cA$-weaved factorization of $f$ is a factorization
\[
X=\lim(Z_j) \to\cdots \to Z_{j+1} \xto{f_{j}} Z_{j} \xto{f_{j-1}}Z_{j-1}\to \cdots\to \mathrm{colim}(Z_j)=Y
\]
with $j$ ranging over the integers, $\cofib(f_j)\in \cA[j]$ for any $j\in \Z $ and with $f_m$ being an isomorphism for $|j|\gg 0$.
\end{definition}

\begin{proposition}
Let $\cA$ be an additive full $\infty$-subcategory of a stable $\infty$-category $\C$ such that any morphism $f\colon X\to Y$  in $\C$ has an $\cA$-weaved factorization. Then the collection of subcategories $\{\cA[n]\}_{n\in \Z}$ is a Bridgeland $\Z$-slicing of $\C$.
\end{proposition}
\begin{proof} Looking at \adef\refbf{def.bridgeland-slicing}, the only thing we need to prove is that $\cA[n_1]\orth\cA[n_2]$ for $n_1>n_2$. This is Remark \refbf{orthogonality-of-a}.
\end{proof}

From Proposition \refbf{b-is-discrete} and Remark \refbf{rem:b-is-discrete} we then immediately have the following converse of Corollary \refbf{cor.oss.for.Heart.to.t}, corresponding to \cite[Lemma 3.2]{Brid}.
\begin{proposition}\label{to.be.repeated.verbatim}
Let $\cA$ be a full additive $\infty$-subcategory of a stable $\infty$-category $\C$, such that  any morphism in $\C$ has a $\cA$-weaved factorization, and let
$\mathcal{H}^n_B\colon \C\to \cA[n]$ be the functors given by taking the cofibers of the $n$-th morphism in the $\cA$-weaved factorization of the initial morphisms. 
 Let
$\C_{\cA,\geq 0}$ be the full subcategory of $\C$ on those objects $X$ such that $\mathcal{H}^j_B(X)=0$ for any $j< 0$, and let $\C_{\cA,< 0}$ be the full subcategory of $\C$ on those objects $X$ such that $\mathcal{H}^j_B(X)=0$ for any $j\geq 0$. Then $\mathfrak{t}_{\cA}=(\C_{\cA,\geq 0}, \C_{\cA,<0})$ is a $t$-structure on $\C$, the stable $\infty$-category $\C$ is bounded with respect to $\mathfrak{t}_{\cA}$, and the standard heart of $\mathfrak{t}_{\cA}$ is (equivalent to) $\cA$. In particular $\cA$ is abelian.
\end{proposition}
\begin{proof}
The only thing left to prove is that the standard heart of  $\mathfrak{t}_\cA$ is $\cA$. To see this notice that an object $Y$ lies in $\C^\heart=\C_{\cA,\geq 0}\cap \C_{\cA,<1}$ if and only if the $\cA$-weaved factorization
\[
0 =\lim(Y_j)\to\cdots \to Y_{j+1} \xrightarrow{f_{j}} Y_{j} \xrightarrow{f_{j-1}}Y_{j-1}\to \cdots\to \mathrm{colim}(Y_j)=Y
\]
of its initial morphism has $\mathrm{cofib}(f_j)=0$ for every $j\neq0$, and so it is of the form
\[
\cdots 0 \to 0\to \cdots \to 0\xrightarrow{f_{0}} Y \xrightarrow{\mathrm{id}_Y}Y\xrightarrow{\mathrm{id}_Y} \cdots\xrightarrow{\mathrm{id}_Y}Y\xrightarrow{\mathrm{id}_Y}\cdots,\] 
with $Y=\mathrm{cofib}(f_0)\in \cA$.
\end{proof}
\section{Semi-orthogonal decompositions.}\label{sec:sods}
\begin{modifyepigraph}{.9}
\epigraph{La vie c'est ce qui se décompose à tout moment; c'est une perte monotone de lumière, une dissolution insipide dans la nuit, sans sceptres, sans auréoles, sans nimbes.}{E. Cioran}
\end{modifyepigraph}

In the previous section we have investigated the case when the equivalence relation $\sim$ from Lemma \refbf{equivalence} had a single equivalence class.
At the opposite end is the case when each equivalence class consists of a single element. As $x\sim x+1$ for any $x\in J$, this is equivalent to requiring that the $\Z$-action is trivial. As noticed in Remark \refbf{rem.finite} this in particular happens when $J$ is a finite finite totally ordered set.
 As we are going to show, this is another well investigated case in the literature: $J$-families of $t$-structures with a finite $J$ capture 
  the notion of \emph{semi-orthogonal decompositions} for the stable $\infty$-category $\C$ (see \cite{Bondal1995, Kuz} for the notion of semi-orthogonal decomposition in the classical triangulated context).

To fix notations for this section, let $J=\ordered{k}$ be the totally ordered set on $(k+1)$ elements, \ie, $J=\{0,1,\dots,k\}$, and let $\tee\colon \mathcal{O}(\ordered{k})\to \ts(\C)$ be a $\ordered{k}$-slicing on $\C$. We have the maximal interval decomposition on $\ordered{k}$ given by intervals $I_i=\{i\}$ for $i=0,\dots,k$. This corresponds to the morphism of posets $\ordered{k-1}\to \mathcal{O}(\ordered{k})$ given by $i\mapsto \{x\geq i+1\}$.  Using this decomposition we have that any morphism $f\colon X\to Y$ in $\C$ has a unique factorization
\[
X \xto{f_{k}} Z_{k-1} \xto{f_{k-1}}Z_{{k-2}}\to\dots\to Z_{{1}} \xto{f_1} Z_{0} \xto{f_0} Y,
\]
with $\cofib(f_i)\in \C_i$, and $\C_i\orth \C_h$, for any $0\leq h <i\leq k$. 

What we are left to investigate are therefore the special features of the $t$-structures $\tee_{i}=(\C_{< i},\C_{\geq i})$ coming from the triviality of the $\Z$-action on $\ordered{k}$, and so on $\mathcal{O}(\ordered{k})$.
By $\Z $-equivariancy of the map $\tee\colon \mathcal{O}(\ordered{k})\to \ts(\C)$, this implies that all the $t$-structures $\tee_{i}$ are $\Z $-fixed points for the natural $\Z $-action on $ \ts(\C)$. Now, a rather pleasant fact is that fixed points of the $\Z$-action on $\ts(\C)$ are precisely those $t$-structures $(\cate{L}, \cate{U})$ for which $\cate{U}$ is a stable sub-$\infty$-category of $\C$. We will make use of the following
\begin{lemma}\label{magicstable}
Let $\cB$ be a full sub-$\infty$-category of the stable $\infty$-category $\C$; then, $\cB$ is a stable sub-$\infty$-category of $\C$ if and only if $\cB$ is closed under shifts in both directions and under pushouts in $\C$.
\end{lemma}
\begin{proof}
The ``only if'' part is trivial, so let us prove the ``if'' part.
First of all let us see that under these assumptions $\cB$ is closed under taking fibers of morphisms. This is immediate: if $f\colon X\to Y$ is an arrow in $\cB$ (\ie an arrow of $\C$ between objects of $\cB$, by fullness), then $f[-1]$ is again in $\cB$ since $\cB$ is closed with respect to the left shift. Since $\cB$ is closed under pushouts in $\C$, also  $\fib(f)=\cofib(f[-1])$ is in $\cB$. It remains to show how this implies that $\cB$ is actually stable, \ie it is closed under all finite limits and satisfies the pullout axiom. Unwinding the assumptions on $\cB$, this boils down to showing that in the square
\[
\xymatrix{
B \ar[r]\ar[d] \pb &  X\ar[d]^f \\
Y \ar[r]_g& Z
}
\]
the pullback $B$ of $f,g \in \hom(\cB)$ done in $\C$ is actually an object of $\cB$; indeed, once this is shown, the square above will satisfy the pullout axiom in $\C$, 
so \emph{a fortiori} it will have the universal property of a pushout in $\cB$. To this aim, let us consider the enlarged diagram of pullout squares in $\C$
\[
\xymatrix{
Z[-1]\ar[r]\ar[d] \ar@{}[dr]|\star & \fib(g)\ar[r]\ar[d] & 0\ar[d] \\
\fib(f)\ar[r]\ar[d] & B\ar[r]\ar[d] & X \ar[d]^f\\
0\ar[r] & Y \ar[r]_g & Z.
}
\]
The objects $Z[-1], \fib(f)$ and $\fib(g)$ lie in $\cB$ by the first part of the proof, so the square $(\star)$ is in particular a pushout of morphism in $\cB$; by assumption, this entails that $B\in\cB$.
\end{proof}
\begin{remark}\label{oss.shifts.pullback}
Obviously, a completely dual statement can be proved in a completely dual fashion:  a full sub-$\infty$-category $\cB$ of a stable $\infty$-category $\C$ is a stable sub-$\infty$-category if and only if it is closed under shifts in both directions and under pullbacks in $\C$.
\end{remark}
\begin{proposition}\label{stableare}
Let $\tee=(\cate{L},\cate{U})$ be a $t$-structure on a stable $\infty$-category $\C$, and let $(\E,\M)$ be the normal torsion theory associated to $\tee$; then the following conditions are equivalent:
\begin{enumerate}
\item $\tee$ is a fixed point for the $\Z $-action on $\ts(\C)$, \ie, $\tee[1]=\tee$ (or equivalently, $\cate{L}[1]= \cate{L}$, or equivalently $\cate{U}[1]=\cate{U}$);
\item $\cate{U}$ is a stable sub-$\infty$-category of $\C$.
\item $\cate{L}$ is a stable sub-$\infty$-category of $\C$.
\item $\E$ is closed under pullback;
\item $\M$ is closed under pushout.
\end{enumerate}
\end{proposition}
\begin{proof}
`(2) implies (1)' is obvious. Namely, if  $\cate{U}$ is a stable sub-$\infty$-category of $\C$, then it is closed under shifts in both directions. Therefore $\cate{U}[-1]\subseteq \cate{U}$, and so $\cate{U}\subseteq \cate{U}[1]$. Since, by definition of $t$-structure, $\cate{U}[1]\subseteq \cate{U}$, we have $\cate{U}[1]= \cate{U}$. To prove that `(1) implies (2)', assume $\cate{U}[1]=\cate{U}$. This implies $\cate{U}[-1]= \cate{U}$, so $\cate{U}$ is closed under shifts in both directions. By Lemma \refbf{magicstable},  we then have only to show that $\cate{U}$ is closed under pushouts in $\C$ to conclude that $\cate{U}$ is a stable $\infty$-subcategory of $\C$. 
Consider a pushout diagram
\[
\xymatrix{
 A \ar[r]\ar[d] \po & B\ar[d] \\
 C \ar[r] & P
}
\]
in $\C$ with $A$, $B$ and $C$ in $\cate{U}$. Since $A$ and $C$ are in $\cate{U} = 0/\E$ we have that both $0\to A$ and $0\to C$ are in $\E$. But $\E$ has the 3-for-2 property, so also $A\to C$ is $\E$. Since $\E$ is closed for pushouts, this implies that also $B\to P$ is in $\E$. But $0\to B$ in in $\E$ since $B$ is in $\cate{U}$, and therefore also $0\to P$ is in $\E$, \ie, $P$ is in $\cate{U}$. 

We now prove that `(1) is equivalent to (4)'. Assume $\E$ is closed under pullbacks. Then for any $X$ in $\cate{U}$ we have that $0\to X$ is in $\E$, and so $X[-1]\to 0$ is in $\E$. By the Sator lemma this implies that $0\to X[-1]$ is in $\E$, \ie, that $X[-1]$ is in $\cate{U}$. This shows that $\cate{U}[-1]\subseteq \cate{U}$ and therefore that $\tee[1]= \tee$. Conversely, assume $\tee[1]=\tee$,
and consider a morphism $f\colon X\to Y$ in $\E$. For any morphism $B\to Y$ in $\C$
consider the diagram
\[
\xymatrix{
\mathrm{fib}(f)\ar[r]\ar[d] & A\ar[r]\ar[d] & X\ar[r]\ar[d]^{f} & 0\ar[d]\\
0 \ar[r] & B\ar[r]  & Y\ar[r]  & \cofib(f)
}
\]
where all the squares are pullouts in $\C$. Since $f$ is in $\E$ and $\E$ is closed for pushouts, also $0\to \cofib(f)$ is in $\E$. This means that $\cofib(f)$ is in $\cate{U}$ and so, since we are assuming that $\cate{U}=\cate{U}[-1]$, also $\mathrm{fib}(f)=\cofib(f)[-1]$ is in $\cate{U}$, \ie,  $0\to \mathrm{fib}(f)$ is in $\E$. By the Sator lemma, $\mathrm{fib}(f)\to 0$ is in $\E$, which is closed for pushouts, and so $A\to B$ is in $\E$. The proofs that `(1) if and only if (3)' and `(1) if and only if (5)' are perfectly dual.
\end{proof}

\begin{remark}\label{oss.hereditary}
A factorization system $(\E,\M)$ for which the class $\E$ is closed for pullbacks is sometimes called an \emph{exact reflective} factorization, see, \eg, \cite{CHK}. This is equivalent to saying that the associated reflection functor is left exact (this is called a \emph{localization} in the jargon of \cite{CHK}). Dually,  one characterizes \emph{co-localizations} of a category $\C$ with an initial object as \emph{co-exact coreflective} factorizations where the right class $\M$ is closed under pushouts.  Therefore, in the stable $\infty$-case, we see that a $\Z $-fixed point in $\ts(\C)$ is a $t$-structure $(\cate{L},\cate{U})$ such that the truncation functors $S\colon \C\to \cate{U}$ and $R\colon \C\to \cate{L}$ respectively form a co-localizations and a localization of $\C$. In the terminology of \cite{Beligiannisreiten} we therefore find that in the stable $\infty$-case $\Z $-fixed point in $\ts(\C)$ correspond to \emph{hereditary torsion pairs} on $\C$. Since we have seen that for a $\Z $-fixed point in $\ts(\C)$ both $\cate{L}$ and $\cate{U}$ are stable $\infty$-categories, this result could be deduced also from \cite[Prop. \textbf{1.1.4.1}]{LurieHA}: a left (resp., right) exact functor between stable $\infty$-categories is also right (resp., left) exact.
 \end{remark}

We can now precisely relate semi-orthogonal decompositions in a stable $\infty$-category $\C$ to $\ordered{k}$-slicings of $\C$. The only thing we still need is the following definition, which is an immediate adaptation to the stable setting of the classical definition given for triangulated categories (see, \eg, \cite{Bondal1995, Kuz} ).
\begin{definition}
Let $\C$ be a stable $\infty$-category. A \emph{semi-orthogonal decomposition} with $k+1$ classes on $\C$ is the datum of $k+1$ stable $\infty$-subcategories $\C_0$, $\C_1$,\dots, $\C_{k}$ of $\C$ such that
\begin{enumerate}
\item one has $\C_i\orth \C_h$ for $h<i$ (semi-orthogonality);
\item for any object $Y$ in $\C$ there exists a unique $\{\C_i\}$-weaved tower, \ie, a factorization of the initial morphism $0\to Y$ as 
\[
0 \xto{f_{k}} Y_{k-1} \xto{f_{k-1}}Y_{{k-2}}\to\dots\to Z_{{1}} \xto{f_1} Y_{0} \xto{f_0} Y,
\]
with $\cofib(f_i)\in \C_i$ for any $i=0,\dots, k$. 
\end{enumerate} 
\end{definition}
Since $\{\C_i\}$-weaved towers are preserved by pullouts, one can equivalently require that any morphism $f\colon X\to Y$ in $\C$ has a unique factorization of has a unique factorization
\[
X \xto{f_{k}} Z_{k-1} \xto{f_{k-1}}Z_{{k-2}}\to\dots\to Z_{{1}} \xto{f_1} Z_{0} \xto{f_0} Y,
\]
with $\cofib(f_i)\in \C_i$, and this immediately leads to the following 
\begin{proposition}\label{what.s.semiortho}
Let $\C$ be a stable $\infty$-category. Then the datum of a semi-orthogonal decompositions with $k+1$ classes on $\C$ is equivalent to the datum of a $\ordered{k}$-slicing on $\C$.
\end{proposition}
\begin{proof}
The only missing piece of information to show that a $\ordered{k}$-slicing is a semi-orthogonal decompositions is the fact that the sub-$\infty$-categories $\C_i$ are stable. But $\C_i=\cate{L}_{i+1}\cap \cate{U}_i$ and both $\cate{L}_{i+1}$ and $\cate{U}_i$ are stable by \aprop\refbf{stableare}. Therefore, also $\C_i$ is stable (see \cite{LurieHA}).  Conversely, given a semi-orthogonal decomposition this defines a $\ordered{k}$-slicing by means of the cofiber functors $\mathcal{H}^i_B\colon \C\to \C_i$, by the same argument in the proof of \aprop\refbf{to.be.repeated.verbatim}. 
\end{proof}
\begin{remark}By Remark \refbf{oss.hereditary}, we recover in the stable $\infty$-setting the well known fact (see \cite[\textbf{IV.4}]{Beligiannisreiten}) that semi-orthogonal decompositions with a single class correspond to \emph{hereditary torsion pairs} on the category.
\end{remark}
\aprop\refbf{what.s.semiortho} immediately suggests to generalize the definition of semi-orthogonal decomposition to the case of an arbitrary toset of indices, not necessarily finite.
\begin{definition}
Let $I$ be a toset, and let $I^\flat$ be the $\Z$-toset given by $I$ endowed with the trivial $\Z$-action (see \refbf{adjoint}). An $I^\flat$-slicing of a stable $\infty$-category $\C$  is called a $I$-semi-orthogonal decomposition of $\C$.
The class of all $I$-semi-orthogonal decompositions of $\C$ will be denoted by $I\textsc{-sod}(\C)$, i.e. $I\textsc{-sod}(\C)=\slicings(I^\flat,\C)$. 
\end{definition}
\begin{remark}
If an $I$-semi-orthogonal decomposition of $\C$ is given, then all the subcategories $\C_i$ are stable, for any $i$ in $I$.
\end{remark}
\begin{remark}\label{J.induces.I}
Let $J$ be a $\mathbb{Z}$-toset, and let $\iota(J)$ be the toset of equivalence classes of $J$, for the equivalence relation $\sim$ of Lemma \refbf{equivalence}. Then every $J$-slicing of a stable $\infty$-category $\C$ induces an $\iota(J)$-semi-orthogonal decomposition of $\C$. Namely, by  \aprop\refbf{adjoint}, $J\rightsquigarrow \iota(J)$ is the left adjoint of the fully faithful embedding $(\firstblank)^\flat\colon \Tos\to \Z\text{-}\Tos$, and the  the projection to the quotient is a $\Z$-equivariant morphism \[
J\to \iota(J)^\flat
\]
which is the unit of this adjunction. By functoriality of the slicings (Remark \refbf{rem.slicing-functor}) we therefore have a natural map
\[
\slicings(J,\C)\to \iota(J)\textsc{-sod}(\C).
\]
\end{remark}

\newcommand{\tilted}{\mathrel{\lightning}\!}
\section{Abelian slicings and tiltings}\label{tiltings}
\begin{modifyepigraph}{.9}
\epigraph{Quando si vuole uccidere un uomo bisogna colpirlo al cuore, e un Winchester è l'arma più adatta.}{R\@. Rojo}
\end{modifyepigraph}
We now review the abelian counterpart of the notion of $J$-slicing and relate slicings on hearts of a stable $\infty$-category $\C$ with slicings of $\C$. First of all recall the notion of \textit{torsion pair} on an abelian $\infty$-category, which is the abelian counterpart of the notion of $t$-structure on a stable $\infty$-category.
\begin{definition}[torsion theory on an abelian $\infty$-category]
Let $\cA$ be an abelian $\infty$-category.
A torsion pair on an abelian $\infty$-category $\cA$ is a pair $(\cate{F},\cate{T})$ of full sub-$\infty$-subcategories of $\cA$ satisfying:
 \begin{enumerate}[label=$\roman*$)]
\item orthogonality: $\cA(X,Y)$ is contractible for each $X \in \cate{T}$, $Y \in \cate{F}$;
\item any object $X \in \cA$ fits into a pullout diagram
\[
\xymatrix{
X_{\cate{T}}\ar[d]\ar[r]& X\ar[d]\\
\zero\ar[r]& X_{\cate{F}}
}
\]
with $ X_{\cate{T}} \in \cate{T}$ and $X_{\cate{F}} \in \cate{F}$. 
 \end{enumerate}
 The subcategories $\cate{T}$ and $\cate{F}$ are called the torsion class and the torsion free class, respectively.
  \end{definition}
  \begin{notat}
We denote by $\tot(\cA)$ the set of torsion theories on $\cA$; this set has a natural choice for a partial order: $(\cate{F}_1, \cate{T}_1) \leq (\cate{F}_2, \cate{T}_2)$ if and only if $\cate{T}_2 \subseteq \cate{T}_1$, or equivalently $\cate{F}_1\subseteq \cate{F}_2$.
\end{notat}
The poset $\tot(\cA)$ has a top and a bottom element, given by $(\cA,\zero)$ and $(\zero,\cA)$, respectively. The following definition is directly inspired by \cite{Rud}. 
 \begin{definition}[abelian slicing]
Let $(I, \leq)$ be a poset. An \textit{abelian $I$-slicing} on $\cA$ is a morphism of posets $\T \colon \mathcal{O}(I) \to \tot (\cA)$ that preserve the top and bottom element. The image of $(\Lambda,\Upsilon)\in \mathcal{O}(I)$ by $\T$ will be denoted $(\cA_{\Lambda}, \cA_{\Upsilon})$
 \end{definition}
 \begin{remark}
Notice that, since there is no choice of a shift functor in an abelian $\infty$-category,  there is no $\Z$-action on $I$ or $\Z$-equivariancy condition involved in the above definition.
 \end{remark}

\begin{remark}[The abelian slicings functor]
By analogy with Remark \refbf{rem.slicing-functor}, for any $\infty$-category $\A$ we have a functor $\slicings_\text{ab} \colon \Pos \to \Pos $  mapping a poset $I$ to the poset of abelian $I$-slicings of $\cA$.
\end{remark}

 \begin{lemma}\label{to.get.slicings.on.heart}
   Let  $\tee_0 = (\cate{L}_0, \cate{U}_0)$ be a $t$-structure on $\C$ with heart $\C^\heart$, and let $\tee_1=\tee_0[1]$. Also let $\tee=(\cate{L}, \cate{U})$ be another $t$-structure with $\tee_0 \leq \tee \leq \tee_1$. Then $$\T=(\cate{F},\cate{T}):=(\cate{L} \cap \C^\heart, \cate{U}\cap \C^\heart)$$
   is a torsion theory on $\C^\heart$. 
 \end{lemma}
\begin{proof}
Clearly, $\cate{F}\subseteq\C^\heart$ and $\cate{T}\subseteq \C^\heart$.
Moreover, $\cate{F}\subseteq\cate{L}$ and $\cate{T}\subseteq \cate{U}$, and so $\cate{T}\orth\cate{F}$.
Now, pick $X \in \C^\heart$ and consider the fiber sequence 
\[
\xymatrix{
X_{\cate{U}}\ar[d]\ar[r]& X\ar[d]\\
\zero\ar[r]& X_{\cate{L}}
}
\]
induced by the $t$-structure $\tee$. From it we get the fiber sequence
\[
\xymatrix{
X_{\cate{L}}[-1]\ar[d]\ar[r]& \zero\ar[d]\\
X_{\cate{U}}\ar[r]& X
}
\]
We have $X_{\cate{L}}[-1]\in \cate{L}[-1]\subseteq \cate{L}\subseteq \cate{L}_1$ and $X\in \C^\heart\subseteq \cate{L}_1$. Since $\cate{L}_1$ is closed by extensions (see Remark \refbf{closed.by.extensions}), this implies that $X_{\cate{U}}\in \cate{L}_1$. Therefore $X_{\cate{U}}\in \cate{L}_1\cap \cate{U}=\cate{L}_1\cap \cate{U}_0\cap \cate{U}=\cate{T}$. An analogous argument shows that $X_{\cate{L}}\in \cate{F}$.
\end{proof}
 \begin{proposition}\label{J-to-t}
Let $(J,\leq)$ be a $\mathbb{Z}$-toset, and let $J^\heart$ a heart of $J$. Then a $J$-slicing on $\C$ induces a $t$-structure on $\C$ together with an abelian $J^\heart$-slicing on $\C^\heart$. 
 \end{proposition}
 \begin{proof}
   Let $ \mathcal{O}(J) \to \ts(\C)$ be a fixed $J$-slicing on $\C$, and let $J^\heart=U_0\cap L_1$ for some (unique) upper set $U_0$ and lower set $L_1$ in $J$. Finally, let $\tee_0$ be the $t$-structure on $\C$ corresponding to the slicing $(L_0,U_0)$ of $J$. Then we know from Remark \refbf{oss.Z.Postnikov} that $\C^\heart= \C_{J^\heart}$ is the standard heart of $\tee_0$. Let $\tee_1$ be the $t$-structure on $\C$ corresponding to the slicing $(L_1,U_1)$ of $J$. By Lemma \refbf{lemma.plus.one} we know that $\tee_1=\tee_0[1]$. Moreover we know from Remark \refbf{oi.vs.oj} that every upper set $\Upsilon$ of $J^\heart$ is of the form $\Upsilon=U\cap J^\heart$ for a unique upper set $U$ in $J$ with $U_0\leq U\leq U_1$. Let $(L,U)$ be the slicing of $J$ determined by $U$. By  Lemma \refbf{to.get.slicings.on.heart},
 \[
 (\Lambda,\Upsilon) \mapsto (\C_L\cap \C^\heart,\C_U\cap \C^\heart) 
 \]
 defines an abelian $J^\heart$-slicing on $\C^\heart$.
\end{proof}
As we are going to show, in the bounded case we also have the converse of the above proposition.

\begin{lemma}\label{verso.il.tilting}
Suppose that $\tee$ is a bounded $t$-structure on $\C$ with heart $\C^\heart$. Then a torsion theory $\T=(\cate{F},\cate{T})$ on $\C^\heart$ induces a bounded $\mathbb{Z} \times_{\mathrm{lex}} \ordered{1}$-slicing on $\C$. 
\end{lemma} 
\begin{proof}
Since every interval of the form $[n_0,n_1]$ in $\mathbb{Z} \times_{\mathrm{lex}} \ordered{1}$ is finite, a bounded $\mathbb{Z} \times_{\mathrm{lex}} \ordered{1}$-slicing is discrete of finite type. Therefore, by \aprop\refbf{b-is-discrete} we are reduced to showing that a torsion theory $\T$ on $\C^\heart$ induces a Bridgeland $\mathbb{Z} \times_{\mathrm{lex}} \ordered{1}$-slicing on $\C$. Since $\T=(\cate{F}, \cate{T})$ is a torsion theory of $\C^\heart$, we have that  $\T[n]=(\cate{F}[n], \cate{T}[n])$ is a torsion theory of $\C^\heart[n]$ for any $n\in \mathbb{Z}$.
Consider the full subcategories
\[
\C_{(n,0)}=\cate{F}[n];\qquad 
\C_{(n,1)}=\cate{T}[n].
\]
Since the $\Z$-action on $\Z\times_{\mathrm{lex}}\ordered{1}$ is $(n,i)+1=(n+1,i)$, it is immediate to see that $\C_{(n,i)+1}=\C_{(n,i)}[1]$ for every $(n,i)$ in $\Z\times_{\mathrm{lex}}\ordered{1}$. Let now $X\in \C_{(n_1,i_1)}$ and $Y\in \C_{(n_2,i_2)}$ with $(n_1,i_1)>(n_2,i_2)$. Since the order is the lexicographic one, we either have $n_1>n_2$ or $n_1=n_2$ and $i_1=1$ and $i_2=0$. In the first case $X\in \C^\heart[n_1]$ and $Y\in \C^\heart[n_2]$ with $n_1>n_2$ and so $X\orth Y$; in the second case $X\in \cate{T}[n_1]$ and $Y\in \cate{F}[n_1]$ and so again $X\orth Y$. Finally, we have to show that for every object $X$ of $\C$ we have a factorization of the initial morphism $0\to X$ into morphisms whose cofibers are in $\C_{n,i}$ for a decreasing sequence of indices $(n,i)$'s in the lexicographic order on $\mathbb{Z} \times \ordered{1}$. To see this, let $X$ be an object in $\C$ and consider the $\C^\heart$-weaved tower of its initial morphism. Keeping only the nontrivial morphisms in this tower we are reduced to a finite factorization of the form
$$0=X_0 \xrightarrow{\alpha_1} \cdots \xrightarrow{\alpha_k} X_k=X$$ 
with $\cofib(\alpha_l)\in \C^\heart[n_l]$ for a suitable sequence of decreasing integers $n_0>n_1>\cdots > n_k$. Since $\T[n_l]=(\cate{F}[n_l], \cate{T}[n_l])$ is a torsion theory on $\C^\heart[n_l]$,  for each $l$ we have a pullout diagram
\[
\xymatrix{
\cofib(\alpha_l)_{\cate{T}[n_l]}\ar[d]\ar[r]& \cofib(\alpha_l)\ar[d]\\
\zero\ar[r]& \cofib(\alpha_l)_{\cate{F}[n_l]}
}
\]
in $\C^\heart[n_l]$. By \aprop\refbf{pullout.is.pullout}, this is a pullout diagram in $\C$ and so we can consider the commutative diagram
\[
\xymatrix{
X_{n_l-1}\ar[r]^{\alpha_{l,+}}\ar[d] & \tilde{X}_{n_l}\ar[d]\ar[r]^{\alpha_{l,-}}&X_{n_l}\ar[d]\\
 \zero\ar[r]&\cofib(\alpha_i)_{\cate{T}[n_l]}\ar[d]\ar[r]& \cofib(\alpha_l)\ar[d]\\
& \zero\ar[r]& \cofib(\alpha_l)_{\cate{F}[n_l]}
}
\]
where each square is a pullout in $\C$. As $l$ ranges from $0$ to $k$  the sequence
\[
  0=X_0 \xrightarrow{\alpha_1,+}\tilde{X}_0 \xrightarrow{\alpha_1,-}X_1\xrightarrow{\alpha_2,+} \cdots \xrightarrow{\alpha_k,-} X_k=X
\]
gives a factorization of $0\to X$ into morphisms whose cofibers are in $\C_{n,i}$ for a decreasing sequence of indices $(n,i)$'s.
\end{proof}
\begin{remark}\label{rem.explicit.tilt}
The nontrivial upper sets in $\Z \times_{\mathrm{lex}} \ordered{1}$ are easily described: they are all of the form $[(n,i),+\infty)$ for some $n\in \Z$ and $i\in \ordered{1}$. In the notation of Lemma \refbf{verso.il.tilting}, the $t$-structures on $\C$ corresponding to these upper sets are easily described by means of Remark \refbf{rem:b-is-discrete}. We have
\begin{align*}
(\C\tilted\T)_{\geq (n,0)} &=\C_{\geq n}=\{X \in \C \mid \mathcal{H}^iX=\zero \text{ for } i<n\}\\
(\C\tilted\T)_{\geq (n,1)} &=\{X \in \C \mid \mathcal{H}^iX=\zero \text{ for } i<n, \; (\mathcal{H}^nX)_{\cate{F}[n]}=\zero\}
\end{align*}
Equivalently,
\[
(\C\tilted\T)_{\geq (n,1)}=\{X \in \C \mid \mathcal{H}^iX=\zero \text{ for } i<n,\; \mathcal{H}^nX\in \cate{T}[n]\}
\] 
\end{remark}

\begin{definition}\label{def.tilting}
Let $\tee=(\C_{<0},\C_{\ge 0})$ be a bounded $t$-structure on $\C$ and let $\C^\heart$ be its standard heart. For every torsion theory $\T$ on $\C^\heart$, the $t$-structure $\tee\tilted\T  = \Big((\C\tilted\T)_{< (0,1)} , (\C\tilted\T)_{\ge (0,1)}\Big)$ is said to be obtained \emph{tilting} (it's a verb) $(\C_{ < 0},\C_{\ge 0})$ with $\T$ (see \cite{Beligiannisreiten}).
\end{definition}

\begin{remark}\label{trivial.tilting}
One immediately sees that the tilting of a $t$-structure $\tee=(\C_{< 0},\C_{ \ge 0})$ by the bottom torsion theory $\T_\perp = (\zero,\C^\heart)$ is the trivial tilting, while tilting by the the top torsion theory $\T_\top = (\C^\heart,\zero)$
correspond to shifting by 1:
\begin{itemize}
 	\item $\tee\tilted\T_\perp = 
 	(\C_{< 0},\C_{\ge 0})$,
 	\item $\tee\tilted\T^\top  = 
 	(\C_{< 1},\C_{\ge 1})=(\C_{< 0},\C_{\ge 0})[1]$.
 \end{itemize} 
\end{remark}

\begin{remark}\label{tilting-as-morphism}
The construction of Lemma \refbf{verso.il.tilting} gives  a map
\[
\tot(\C^\heart)\to \slicings(\mathbb{Z} \times_{\mathrm{lex}} \ordered{1},\C ),
\]
and the explicit description in Remark \refbf{rem.explicit.tilt} show that this is a morphism of posets. Namely, if $\T_1\leq \T_2$ then $\cate{T}_2[n] \subseteq \cate{T}_1[n]$ and so $(\C\tilted\T_2)_{\geq (n,i)} \subseteq (\C\tilted\T_1)_{\geq (n,i)}$ for every $(n,i)$. In particular, tilting defines a morphism of posets
\[
(\secondblank \tilted \firstblank) \colon \tot(\C^\heart)\times \ts(\C) \to \ts(\C).
\]
This construction can be seen as a byproduct of the functoriality of slicings as follows. The $\mathbb{Z}$-toset $\mathbb{Z} \times_{\mathrm{lex}} \ordered{1}$ has an obvious $\mathbb{Z}$-equivariant morphism of posets  to $\mathbb{Z}$ given by the projection $\pi$ on the first factor. However, and remarkably, there is also another less trivial $\mathbb{Z}$-equivariant morphism 
\[
\pi\circ \mathrm{tilt}\colon \mathbb{Z} \times_{\mathrm{lex}} \ordered{1} \to \mathbb{Z}
\]
given by the composition of the projection $\pi$ with the  $\mathbb{Z}$-equivariant automorphism of the poset $\mathbb{Z} \times_{\mathrm{lex}} \ordered{1}$ defined by
\[
\begin{cases} \mathrm{tilt}(n,0)=(n-1,1) \\ \mathrm{tilt}(n,1)=(n,0) \end{cases}. 
\]
Since taking slicings of a fixed stable $\infty$-category with respect to a $\Z$-toset $J$ is functorial in $J$ (Remark \refbf{rem.slicing-functor}), we have a morphism of $\Z$-tosets
\[
\slicings(\mathbb{Z} \times_{\mathrm{lex}} \ordered{1},\C ) \xrightarrow{\pi\circ\mathrm{tilt}} \ts(\C).
\]
Composing this with the morphism of posets $\tot(\C^\heart)\to \slicings(\mathbb{Z} \times_{\mathrm{lex}} \ordered{1},\C )$ gives the tilting map.
\end{remark}

The proposition below can be found in \cite{Beligiannisreiten} for $J=\mathbb{Z} \times_{\rm{lex}} \ordered{1}$ and in \cite{Brid} for $J=\mathbb{R}$. 
\begin{proposition}\label{converse.if.finite}
Let $J$ be a $\mathbb{Z}$-toset, and let $J^\heart$ be a heart of $J$. Then giving a bounded $J$-slicing on $\C$ is equivalent to giving a bounded $t$-structure $(\C_{< 0},\C_{\ge 0})$ on $\C$ together with an abelian $J^\heart$-slicing on the standard heart $\C^\heart$ of $(\C_{< 0},\C_{\ge 0})$. Moreover the bounded $J$-slicing on $\C$ is discrete if and only if the abelian $J^\heart$-slicing of $\C^\heart$ is discrete.
\end{proposition}
\begin{proof}
In one direction this is the content of \aprop\refbf{J-to-t}. Vice versa, assume we have a bounded $t$-structure $(\C_{< 0},\C_{\ge 0})$ on $\C$ and an abelian $J^\heart$-slicing $\T$ on its standard heart $\C^\heart$. By definition this is a morphism of posets $\mathcal{O}(J^\heart)\to \tot(\C^\heart)$. By Remark \refbf{tilting-as-morphism}, tilting a fixed $t$-structure gives a morphism of posets $\tot(\C^\heart)\to \ts(\C)$, and so by composition we get a morphism of posets 
\[
\mathcal{O}(J^\heart)\to\ts(\C)
\]
Recalling the identification of $\mathcal{O}(J^\heart)$ with the interval $[U_0,U_1]$ of $\mathcal{O}(J)$ from Remark \refbf{oi.vs.oj}, and that $U_1=U_0+1$ from Lemma \refbf{lemma.plus.one}, this is a morphism of posets $[U_0,U_0+1]\to \ts(\C)$ and so it induces a uniquely determined $\Z$-equivariant morphism of $\Z$-tosets
\[
\tee\colon \Z\times_{\mathrm{lex}}[U_0,U_0+1]\to \ts(\C)
\]
By Remark \refbf{trivial.tilting}, $\tee_{(1,U_0)}=\tee_{(0,U_0)}[1]=\tee_{(0,U_0+1)}$ and so $\tee$ factors through the natural morphism of $\Z$-tosets 
\[
\Z\times_{\mathrm{lex}}[U_0,U_0+1]\to \mathcal{O}(J)
\]
given by $(n,U)\mapsto U+n$. In other words, $\tee$ uniquely defines a $J$-slicing on $\C$, 
which is bounded since the $t$-structure $(\C_{< 0},\C_{\ge 0})$ is, by Remark \refbf{bounded.is.bounded}. Finally, the construction manifestly preserves finite types.
\end{proof}
\section{Concluding remarks}
\label{concluding}
\begin{modifyepigraph}{.2}
\epigraph{That's all, folks!}{Bosko}
\end{modifyepigraph}

We have explored two classes of $J$-slicings so far: those for which $J$ has a heart, and those for which $\mathbb{Z}$ acts trivially on $J$. In this section, we show how these two cases are fundamental building blocks for all other $J$-slicings. 
\begin{lemma}\label{for-the-main-theorem}
Let $\tee$ be a $J$-slicing on a stable $\infty$-category $\C$, and let $I_i\subseteq J$ be the equivalence class of $i\in J$ with respect to the equivalence relation $\sim$ of Lemma \refbf{equivalence}. For every $(\Lambda,\Upsilon)=(L\cap I_{i},U\cap I_{i})$ in $\mathcal{O}(I_i)$, let $\tee_{i;\Lambda,\Upsilon}=(\C_L\cap \C_{I_i},\C_U\cap \C_{I_i})$. Then $\tee_i\colon (\Lambda,\Upsilon)\to \tee_{i;\Lambda,\Upsilon}$ is a $I_i$-slicing of $\C_{I_i}$.
\end{lemma}
\begin{proof}
By Remark \refbf{J.induces.I}, the $J$-slicing $\tee$ induces an $\iota(J)$-semi-orthogonal decomposition of $\C$: for every equivalence class $[i]$ in $\iota(J)$ the corresponding slice in this semi-orthogonal decomposition is the subcategory $\C_{I_i}$ of $\C$ determined by the $J$-slicing $\tee$, where $I_i\subseteq J$ is the equivalence class of $i$ with respect to the equivalence reltion $\sim$, as a subset of $J$. As they are the slices of a semi-orthogonal decomposition, the subcategories $\C_{I_i}$ are stable (this can also be seen directly from the definition of the $\C_{I_i}$'s).
 As shown in Example \refbf{class.is.interval}, $I_i$ is an interval of $J$ and a sub-$\Z$-toset of $J$, simply by definition of the equivalence relation. 
For every $i$, we can therefore write $I_i=U_{i;0}\cap L_{i;1}$. By Remark \refbf{oi.vs.oj} every slicing $(\Lambda,\Upsilon)$ of $I_i$ is of the form $\Lambda=L\cap I_{i}$ and $\Upsilon=U\cap I_{i}$ for a unique slicing $(L,U)$ of $J$ with $U_{i;0}\leq U\leq U_{i;1}$. This gives an isomorphism of tosets between $\mathcal{O}(I_i)$ and the interval $[U_{i;0},U_{i;1}]$ in $\mathcal{O}(J)$.
Now, to show that $\tee_{i;\Lambda,\Upsilon}$ is a $t$-structure on $\C_{I_i}$ one verbatim repeats the proof of Lemma \refbf{to.get.slicings.on.heart} to get orthogonality of the classes and the existence of the relevant fiber sequences. Next, to show that $(\C_U\cap \C_{I_i})[1]\subseteq \C_U\cap \C_{I_i}$ notice that, since $I_i$ is an equivalence class, we have $I_i+1=I_i$ and so 
\[
(U\cap I_{i})+1=(U+1)\cap (I_i+1)\subseteq U\cap I_i.
\]
This shows that $\tee_i$ is a morphism of sets $\mathcal{O}(I_i)\to \ts(\C_{I_i})$ and it is immediate to see that this is actually a morphism of tosets. Finally, $\Z$ equivariance is obtained by noticing that $\Upsilon+1=(U\cap I_i)+1=(U+1)\cap I_i$ and $\Lambda+1=(L+1)\cap I_i$, so that
\[
\tee_{i;\Lambda+1,\Upsilon+1}=(\C_{L+1}\cap \C_{I_i},\C_{U+1}\cap \C_{I_i})=(\C_L[1]\cap \C_{I_i},\C_U[1]\cap \C_{I_i})=\tee_{i;\Lambda,\Upsilon}[1],
\]
where we used that $\C_{I_i}$ is a stable subcategory of $\C$ and so $\C_{I_i}=\C_{I_i}[1]$.
\end{proof}

We can now state and prove our main result, summarising and putting together the various pieces constructed so far. To make a self-standing statement, we explicitly recall the definition of the equivalence relation $\sim$ from Lemma \refbf{equivalence} in the statement of the theorem below.
\begin{theorem}\label{conclusion}
Let $(J,\leq)$ a $\Z$-toset. A finite type $J$-slicing on a stable $\infty$-category $\C$ is equivalent to the following data:
\begin{enumerate}
\item[(i)] a finite type semi-orthogonal decomposition of $\C$ whose slices $\C_{[x]}$ are stable  $\infty$-subcategories of $\C$ indexed by equivalence classes in $J$ with respect to the equivalence relation $x\sim y$ if and only if there exist integers $n_1$ and $n_2$ with $x+n_1\leq y\leq x+n_2$;
\item[(ii)] a bounded $t$-structure $\tee_{[x]}$ on $\C_{[x]}$ for each  $[x]$ in $J/_{\!\sim}$;
\item[(iii)] a finite type abelian $[x,x+1)$-slicing on $\C_{[x]}^\heart$ for every $[x]$ in $J/_{\!\sim}$ such that $x$ is not a fixed point of the $\Z$-action on $J$.
\end{enumerate}
\end{theorem}
\begin{proof}
By Lemma \refbf{for-the-main-theorem}, a $J$-slicing on a stable $\infty$-category induces semi-orthogonal decomposition of $\C$ whose slices $\C_{I_i}$ are stable  $\infty$-subcategories of $\C$ indexed by equivalence classes $I_i$ in $J$ with respect to the equivalence relation $\sim$, together with $I_i$-slicings of these subcategories.
If $i$ is a fixed point for the $\Z$-action on $J$, then $I_i=\{i\}$ and $\mathcal{O}(I_i)=\ordered{1}$ so that a $I_j$-slicing is trivial. On the other hand, by Example \refbf{I.has.a.heart}, precisely when $i$ is not a fixed point of the $\Z$-action the interval $I_i$ has a heart $I_i^\heart$ which can be identified with the interval $[i,i+1)$ of $J$. Therefore, by \aprop\refbf{J-to-t} an $I_i$-slicing on $\C_{I_i}$ induces a $t$-structure on $\C_{I_i}$ together with an abelian $[i,i+1)$-slicing on the standard heart $\C_{I_i}^\heart$. Moreover, by \aprop\refbf{converse.if.finite}, in the finite type case an $I_i$-slicing on $\C_{I_i}$ is precisely equivalent to this datum of a bounded $t$-structure on $\C_{I_i}$ with an abelian $I_i^\heart$-slicing on $\C_{I_i}^\heart$. 
\par
Vice versa, given the data (i)-(iii), we want to define a $J$-slicing of $\C$. For every upper set $U$ of $J$, we can write 
\[
U=\bigcup_{[i]\in J/\sim}U\cap I_i
\]
and, by Prop. \refbf{converse.if.finite} again, the data (i)-(iii) define $\infty$-subcategories $\C_{U\cap I_i}\subseteq \C_{[i]}$ of $\C$. Defining $\C_U$ as the extension closed $\infty$-subcategory of $\C$ generated by the $\C_{U\cap I_i}$'s provides a map $\tee\colon \O(J)\to \ts(\C)$ which is immediate to see it is montone and $\mathbb{Z}$-equivariant and so is a $J$-slicing of $\C$. Moreover, this $J$-slicing is of finite type as all the slicings provided by the data are of finite type, and manifestly this way of constructing  finite type $J$-slicings out of
data (i)-(iii) provides an inverse to the construction described in the first part of the proof.  
\end{proof}

\hrulefill
\subsubsection*{Translations of the opening quotes.}
\footnotesize
\begin{enumerate}[label=\S\arabic*.)]
 \setcounter{enumi}{1}
 \item These are the first three verses of \emph{Tao te Ching}'s chapter 63: ``Act without action / work without work / taste without taste.''
 \item This is a quote from \emph{Histoire d'O}: ``her freedom was worse than any chain.''
 \item This is \emph{Genesis 11:7}: ``let us go down, and there confound their language, that they may not understand one another's speech.''
 \setcounter{enumi}{5}
 \item This is a quote from E. Cioran's \emph{Brief history of decay}: ``Life is what decomposes at every moment; it is a monotonous loss of light, an insipid dissolution in the darkness, without scepters, without halos.''
 \item This is a famous line Ramon Rojo says in S. Leone's movie \emph{A fistful of dollars}: ``When you want to kill a man you must shoot for his heart and a Winchester is the best weapon.''
\end{enumerate}
\bibliographystyle{amsalpha}
\bibliography{allofthem}

\end{document}